\documentclass{siamltex}
\usepackage{amsfonts} 
\usepackage[centertags]{amsmath}
\usepackage{amssymb}
\usepackage{booktabs,color}
\usepackage{newlfont}
\usepackage[margin=1in]{geometry}
\usepackage{graphicx}
\usepackage[figurewithin=none,hang,small,md]{caption}
\usepackage[hang,small,md,TABBOTCAP]{subfigure}
\usepackage[ansinew]{inputenc}
\usepackage{colortbl}
\usepackage{pst-all}
\usepackage{wasysym}
\usepackage{mathrsfs}
\usepackage{float}
\usepackage{verbatim}
\usepackage{upgreek}

\providecommand{\U}[1]{\protect\rule{.1in}{.1in}}

\renewcommand{\eqref}[1]{(\ref{#1})}
\usepackage[numbers,sort&compress]{natbib}

\newcommand{\eps}{{\displaystyle \varepsilon}}
\newcommand{\bsub}{\begin{subequations}}
\newcommand{\esub}{\end{subequations}$\!$}
\newcommand{\ds}[0]{\displaystyle}

\newcommand{\bigoh}{\mathcal{O}}

\newcommand{\bx}{\mathbf{x}}

\newcommand{\by}{\mathbf{y}}

\newcommand{\dS}{dS}
\newcommand{\sig}{\nu}
\newcommand{\Om}{\Omega}
\newcommand{\pOm} {\partial \Omega}
\newcommand{\Gama}{\Gamma_a}
\newcommand{\Gamr}{\Gamma_r}
\newcommand{\Gamj}{\Gamma_j}

\newcommand{\A}{{\mathcal{A}}}
\newcommand{\E}{{\mathcal{E}}}
\newcommand{\Np}{N}

\begin{document}

\title{Numerical approximation of diffusive capture rates by planar and spherical
surfaces with absorbing pores}

\date{\today} 
\author{
Andrew J. Bernoff\thanks{{Dept.~of Mathematics, Harvey Mudd College, Claremont, CA, 91711, USA. {\tt ajb@hmc.edu}}} \and 
Alan E. Lindsay\thanks{Dept.~of Applied and Computational Math \& Statistics, University of Notre Dame, Notre Dame, Indiana, 46656, USA. {\tt a.lindsay@nd.edu} }
}

\baselineskip=16pt

\maketitle

\begin{abstract}
In 1977 Berg \& Purcell published a landmark paper entitled \emph{Physics of Chemoreception} which examined how a bacterium can sense a chemical attractant in the fluid surrounding it \cite{bergp}. At small scales the attractant molecules move by Brownian motion and diffusive processes dominate. This example is the archetype of  \emph{diffusive signaling problems} where an agent moves via a random walk until it either strikes or eludes a target.  Berg \& Purcell modeled the target as a sphere with a set of small circular targets (pores) that can capture a diffusing agent. 
They argued that, in the limit of small radii and wide spacing, each pore could be modeled independently as a circular pore on an infinite plane. Using a known exact solution, they showed the capture rate to be proportional to the combined perimeter of the pores. In this paper we study how to improve this approximation by including inter-pore competition effects and verify this result numerically for a finite collection of pores on a plane or a sphere.  Asymptotically we have found the corrections to the Berg-Purcell formula that account for the enhancement of capture due to the curvature of the spherical target and the inhibition of capture due to the spatial interaction of the pores. Numerically we develop a spectral boundary element method for the exterior mixed Neumann-Dirichlet boundary value problem.  Our formulation reduces the problem to a linear integral equation, specifically a Neumann to Dirichlet map, which is supported only on the individual pores. The difficulty is that both the kernel and the flux are singular, a notorious obstacle in such problems. A judicious choice of singular boundary elements allows us to resolve the flux singularity at the edge of the pore. In biological systems there can be thousands of receptors whose radii are $0.1\%$ the radius of the cell.  Our numerics can now resolve this realistic limit with an accuracy of roughly one part in $10^8$.\end{abstract}


\label{firstpage}

\setcounter{equation}{0}
\setcounter{section}{0}



\pagestyle{myheadings}
\markboth{A.~J.~Bernoff, A.~E.~Lindsay}{Capture Problems on the Plane and the Sphere}

\section{Introduction}\label{sec:intro}

Cells are calculating machines which must infer information about chemical concentrations in their environment and make decisions based on these measurements. For example, T-cell receptors on the cell membrane must trigger an immune response when foreign bodies are encountered \cite{Pageon2016,Rossy2012,zwanzig}.  In their seminal paper \cite{bergp}, Berg \& Purcell studied the fundamental biophysical limits on chemical sensing at microscopic scales and demonstrated that cells could have nearly optimal sensing performance, provided their receptors were numerous and distributed over their exterior membrane. While this leading order theory elucidates the underlying principles of chemoreception, it does not account for the reduction in sensing performance from inter-pore competition. A longstanding problem is to resolve this limitation and describe how the number and detailed spatial arrangement of receptors dictate the ability of the cell to sense its surroundings through diffusive contact. 



\begin{figure}
\centering
\subfigure[Plane with finite cluster of absorbing pores.]{ \label{Fig:schematic_a} \includegraphics[width=0.475\textwidth]{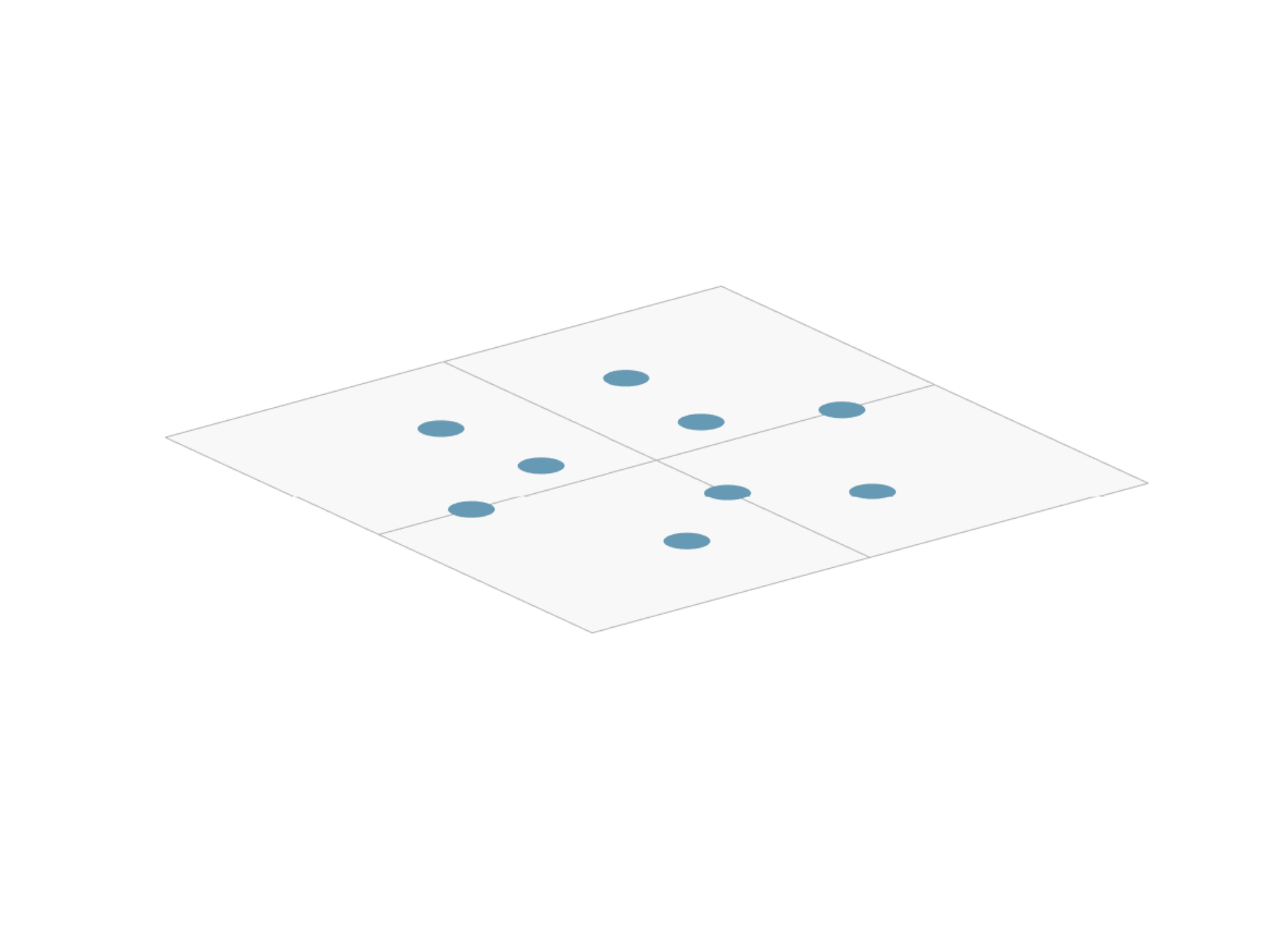}}
\subfigure[Sphere with absorbing surface pores.]{\label{Fig:schematic_b} \includegraphics[width=0.475\textwidth]{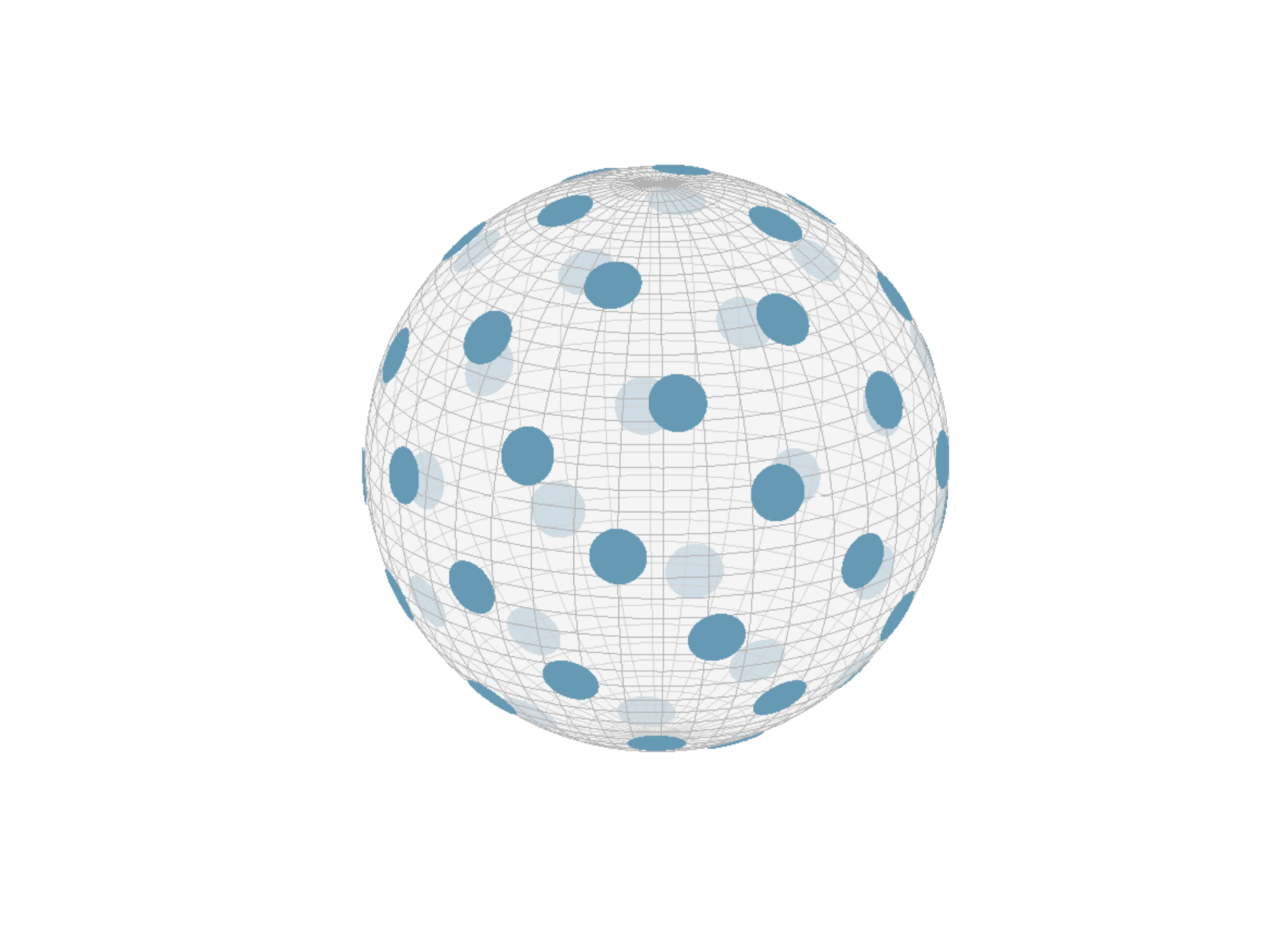}}
\caption{Schematic diagram of the planar and spherical absorption geometries with absorbing surface pores.
\label{Fig:schematic}}
\end{figure}


The mathematical formulation for this problem draws from the theory of electrostatics \cite{sneddon1966mixed,K} and is governed by Laplace's equation with a mixed configuration of Neumann and Dirichlet conditions corresponding to reflecting and absorbing portions of the target. For geometries akin to those in Fig.~\ref{Fig:schematic}, the probability, $v(\bx)$,  that a particle originating at $\bx = (x,y,z)$ in an exterior domain $\Om$ avoids absorption at the target set $\Gamma_a$ on the surface $\pOm$ satisfies
\bsub\label{CapExt}
\begin{gather}
\label{CapExt_a} \Delta v = 0, \qquad  \bx \in \Omega; \\[5pt]
\label{CapExt_b} v = 0, \qquad \bx \in \Gama; \qquad \partial_n v = 0, \qquad \bx \in  \Gamr;\\[5pt]
\label{CapExt_c} v(\bx) =  1 - \frac{C}{|\bx|} + \bigoh\left(\frac{1}{|\bx|^2} \right),  \qquad |\bx| \to \infty.
\end{gather}
\esub
Here the surface $\pOm$ contains an absorbing domain, $\Gama$, which is typically the union of $N$ non-overlapping pores,  $\Gama = \displaystyle\cup_{j=1}^N \Gamj$. The remainder of the surface $\Gamr = \pOm\setminus \Gama$ is reflecting. In the far-field the probability approaches unity as almost all particles will escape capture. The parameter $C$, known as the \emph{capacitance}, is determined uniquely by \eqref{CapExt} and fixes the total flux $J = D\int_{\partial\Omega}\partial_n v\,dS$ of particles to the target.

In the present work, we consider two specific scenarios in which a finite cluster of traps are arranged on an infinite plane (Case I), and where the traps are arranged on the surface of a sphere (Case II) as shown in Fig.~\ref{Fig:schematic}. For both cases we will define the pores as the set of points on the surface $\pOm$ within a distance $\eps a_j$ of some center point $\bx_j$ and seek to understand how the number and spatial arrangement of traps affects the total flux $J$ of particles to the target. The parameter $\eps$ is a common scale factor associated with each pore, and ensures that pores with distinct centers are non-overlapping as $\eps \to 0$.  In Case I, the domain $\Omega$ is the half-space $z>0$ and the absorbing target set is given explicitly as
\begin{equation}\label{PlanarTargets}
 \Gama = \bigcup_{j=1}^N \Gamj,  \qquad \Gamj = \{ (x,y,0) \ | \ (x -x_j)^2 + (y-y_j)^2 < \eps^2 a_j^2 \}
\end{equation}
where $\Gamj$ has center  $\bx_j = (x_j,y_j,0)$.

In Case II, the domain $\Omega$ is the exterior of the unit sphere,  $|\bx| >1$. We use spherical coordinates $(\rho, \theta,\phi)$ on the sphere $\rho=1$ with the pore 
$\Gamj$ centered at $\bx_j = (\sin \theta_j \cos \phi_j, \sin \theta_j \sin \phi_j,\cos \theta_j )$. The absorbing target set is  given explicitly in spherical coordinates as
\begin{equation}\label{SpherePatch}
\Gama =  \bigcup_{j=1}^N \Gamj,  \qquad \Gamj = \{ (1,\theta,\phi ) \ | \ 2\left [ 1-  \sin \theta \,\sin  \theta_j \, \cos(\phi-\phi_j) -\cos\theta \,\cos \theta_j \right ] < \eps^2 a_j^2 \} .
\end{equation} 
An alternative interpretation of $v(\bx)$ is that it is the equilibrium concentration observed in $\Omega$ for a diffusing species when a target is immersed in a uniform concentration that is unity in the far-field. The divergence theorem allows us to relate the flux of the species into the pores to the far-field behavior.
For the planar case, the flux $J_p$ is given by
\begin{equation}\label{FluxFormulaPlane}
J_p = D\int_{\partial\Omega} \partial_n v\, \dS = 2\pi D C,
\end{equation}
while in the spherical case,
\begin{equation}\label{FluxFormulaSphere}
J_s = D\int_{\partial\Omega} \partial_n v\, \dS = 4\pi D C.
\end{equation}

Before stating the main results of this paper, we review some key results associated with these classical problems. The seminal study of Berg and Purcell \cite{bergp} analyzed a spherical target of radius $R_0$ (Case II) partially covered by localized absorbing receptors. From a flux based analysis, they postulated that $N$ non-overlapping receptors of common radius $a_0$, would give rise to a capacitance $C_{\textrm{bp}}$ and associated flux $J_{\textrm{bp}}$  
\begin{equation}\label{resultBP}
J_{\textrm{bp}} = 4\pi D C_{\textrm{bp}}, \qquad C_{\textrm{bp}} = \frac{N a_0 R_0}{N a_0 + \pi R_0}. 
\end{equation}
A key insight from \eqref{resultBP} is that when the pores are well separated on the sphere (which implies $Na_0 \ll R_0$), the flux is proportional to the perimeter of the absorbing set $\Gama$, $J \approx 4 D N a_0$.  Therefore, for fixed absorbing area (which could be a small fraction of the total surface area), a distributed and fragmented absorbing set (with a large perimeter) can have a capture rate that approaches that of an all absorbing target ($J=4 \pi R_0D$).  The expression \eqref{resultBP} was derived through physical reasoning and interpolates between these two limits. The result \eqref{resultBP} does not, however, inform on how the particular spatial arrangement of pores contributes to the capture rate. The microscopic patterning or \emph{clustering} of receptor sites on membrane surfaces is frequently observed experimentally and known to play a key biophysical role in many systems \cite{Care2011,Roh2015,Rossy2012,Pageon2016}. 

Evaluating the impact of clustering on the capture rate is challenging on account of several factors. Exact solutions to \eqref{CapExt} are unfortunately not available beyond the most rudimentary scenarios. Numerical studies of \eqref{CapExt} are also challenging due to the heterogeneous array of mixed Neumann and Dirichlet boundary conditions. Such problems are notorious in potential theory due to a flux discontinuity along the perimeters of the absorbing pores \cite{K,sneddon1966mixed,Duffy08}. Brownian particle simulations \cite{berez2012,berez2013,berez2006,berez2004,Care2011} are a widely used  and flexible approach to sampling the flux $J$ which avoid the challenges of resolving the discontinuous potential, however, these methods are slow to converge, offer relatively crude accuracies and have difficulties dealing with the infinite computational domains and the small pore sizes inherent in biologically realistic applications. 

A complementary approach to tackling the complications of the heterogeneous boundary conditions is to seek a homogenized or effective medium theory. At a large distance from a target surface with an array of pores, the solutions of \eqref{CapExt}  are largely independent to variations in lateral directions. Therefore, the complex configuration of mixed Neumann and Dirichlet conditions can be replaced by a uniform Robin boundary condition $\partial_n v + \kappa\, v = 0$ on the target surface \cite{muratov,shoup,berez2004,berez2006,zwanzig}. Here $\kappa$ is sometimes referred to as the \emph{leakage parameter} as it governs the flux leaking through the boundary. These theories estimate the dependence of the leakage parameter on the absorbing area fraction, $\sigma$, of the surface.
The first boundary homogenization \cite{shoup} was for a sphere (akin to Fig.~\ref{Fig:schematic_b})  based on the Berg-Purcell formula \eqref{resultBP} and specified 
\begin{equation}\label{KShoup}
\kappa = \kappa_{\textrm{bp}} = \frac{4D}{\pi a} \sigma,
\end{equation}
so that the total flux to the target of the homogenized and full problems \eqref{CapExt} were equivalent. Our recent work \cite{LWB2017} reviews the history of homogenization for the sphere and extends this result to incorporate the arrangement and interaction of the pores.

Other studies have considered periodic arrays of traps on a plane bounding a half-space \cite{berez2004,berez2006,muratov}. They have proposed the functional form
\begin{equation}\label{KBerez}
\kappa_{\textrm{be}} = \frac{4 D \sigma}{\pi a}f(\sigma), \qquad f(\sigma) = \frac{1+\alpha\sqrt{\sigma} - \beta \sigma^2}{(1-\sigma)^2},
\end{equation}
where particle simulations estimated the parameter values  for clusters of absorbing pores in square and hexagonal lattices respectively. 

Finite cluster of pores (as depicted in Fig.~\ref{Fig:schematic_a}) have been considered by a set of recent studies \cite{berez2012,berez2013,berez2014} which propose that the cluster be replaced by a single circular pore on which the Robin boundary condition applies. The leakage parameter is again estimated by the formula \eqref{KBerez} where 
$\sigma$ is replaced by the effective pore density within the occupied cluster sites.

The contribution of this work is two-fold. First, we present matched asymptotics formulae for the capacitance of finite sets of pores on the surface of a sphere or on a half-space in the biologically relevant limit of large separation. Second, we present a spectral boundary element that is capable of verifying these results to a high accuracy.   In \S \ref{sec:planarpore}, we derive by matched asymptotic expansions, the following asymptotic expression for the flux $J_p = 2\pi D C$ in the planar case. When all the pores have common radius $\eps$, we find that
\begin{equation}\label{mainResult}
J_p = 4\eps DN \left[ 1- \frac{2\eps}{N\pi} \sum_{j\neq k} \frac{1}{|\bx_j - \bx_k|} + \frac{4\eps^2}{N\pi^2} \sum_{j\neq k}\sum_{i\neq j} \frac{1}{|\bx_j-\bx_k| |\bx_i - \bx_j|} + \bigoh(\eps^3)\right].
\end{equation}
The analogous result for the unit sphere, obtained in \cite{LWB2017}, gives the limiting form of $J_s = 4\pi D C$ as
\bsub\label{FluxSphere}
\begin{equation}\label{FluxSphere_a}
J_s = 4\eps D N\left[ 1 -  
 \frac{\eps}{\pi}\log 2\eps +
 \frac{\eps}{\pi} \Big(\frac{3}{2} -
  \frac{2}{N}\sum_{k\neq j} g_s(|\bx_j-\bx_k|) \Big) + \bigoh(\eps^2\log\eps) \right] \qquad \eps\to0,
\end{equation}
where the spherical pore competition kernel $g_s(\mu)$ is given by
\begin{equation}\label{FluxSphere_b}
g_s(\mu) = \frac{1}{\mu} + \frac{1}{2} \log\Big( \frac{\mu}{2+\mu} \Big), \qquad 0<\mu\leq2.
\end{equation}
\esub
The importance of the asymptotic formulas \eqref{mainResult} and \eqref{FluxSphere} is that they give a first principles account of how the spatial configuration of surface pores affects the capture rate of the target.  The leading order term $4\eps ND$, which appears in both formulae is the classic Berg-Purcell term \eqref{resultBP}, informs on how the perimeter of the pore set influences absorption. The subsequent terms gives corrections due to inter-pore competition, and in the spherical case, logarithmic terms which account for the curvature of the target. 

The rest of the paper is organized as follows. In \S\ref{sec:planarpore}, we use a matched asymptotic expansion analysis to obtain expression \eqref{mainResult}, the flux of diffusing particles to $N$ well separated absorbing pores arranged on a plane. In \S\ref{sec:numerics} we derive a spectral boundary element method for the efficient numerical solution of \eqref{CapExt} which incorporates the known form of the flux singularity and allows for rapid and accurate evaluation of $J$. In \S\ref{sec:results}, we use this method to explore the capture rate of a variety of target sets $\Omega$ ranging up to thousands of pores. Finally, in \S\ref{sec:discussion}, we discuss the implications of the present work and highlight avenues for future investigations.

\section{Asymptotic analysis of the planar pore capture problem}\label{sec:planarpore}

In this section we detail a singular perturbation analysis which yields an approximation for the capacitance of a plane with absorbing circular pores \eqref{mainResult}. One of the first steps is to rephrase the solution of \eqref{CapExt} as $v(\bx) = -C_p u(\bx)$ where $u(\bx)$ satisfies the associated problem 
\bsub\label{CapExtU}
\begin{gather} 
\label{CapExtU_a} \Delta u = 0, \qquad \bx\in \Omega; \\[5pt]
\label{CapExtU_b} u = 0, \qquad \bx \in \Gamma_a, \qquad \partial_n u = 0 \qquad \bx \in \Gamma_r;\\[5pt]
\label{CapExtU_c} u(\bx) =  -\frac{1}{C_p} + \frac{1}{|\bx|}   + \bigoh\left( \frac{1}{|\bx|^2}\right), \qquad |\bx| \to \infty.
\end{gather}
\esub
In the analysis of this section, the domain $\Omega$ and its boundary $\partial\Omega$ are defined as
$$
\Omega = \{ (x,y,z) \in\mathbb{R}^3 \ | \ z >0 \}, \qquad  \partial\Omega = \{ (x,y,z) \in\mathbb{R}^3 \ | \ z =0 \},
$$
and $\partial\Omega = \Gamma_a \cup \Gamma_r$. The formulations \eqref{CapExt} and \eqref{CapExtU} differ in their normalization. In \eqref{CapExt}, $\lim_{|\bx|\to\infty} v(\bx) = 1$  which uniquely determines the capacitance $C_p$ from the strength of the monopole as $|\bx|\to\infty$. In \eqref{CapExtU}, the strength of the monopole in the far field is normalized to unity, implying that $\int_{\partial\Omega}\partial_n u \, \dS = -2\pi$. For the normalization condition \eqref{CapExtU_c}, the capacitance is determined from $\lim_{|\bx|\to\infty}u(\bx) = -C_p^{-1}$.

 Exact solutions to \eqref{CapExtU} have been developed in the simple cases where the absorbing set $\Gamma_a$ is one \cite{sneddon1966mixed} or two non-overlapping absorbing pores \cite{Strieder08,Strieder12}. However, these methods rely heavily on exploiting  symmetries of the set $\Gamma_a$ and cannot be easily generalized for larger $N$. The alternative approach taken here is based of a matched asymptotic analysis \cite{LTS2016,LWB2017,Coombs2009,Ward2010,ChevWard2010} in which we assume the presence of $N$ non-overlapping and well separated pores centered at $\bx_k = (x_k,y_k,0)$ with radii $\eps a_k $ as $\eps\to0$. The parameters $a_k$ allow the pores to have variable radii. A key constituent in the solution of \eqref{CapExtU} is knowledge of the Green's function $G_p(\bx,\bx_0)$ where $\bx_0 = (x_0, y_0,0)$ satisfying
\bsub\label{SurfaceGreens}
\begin{gather}\label{SurfaceGreens_a}
\Delta G_p = 0 ,\quad \bx \in \Omega; \qquad \partial_n G_p = - \delta (\bx-\bx_0), \quad \bx\in\partial\Omega.
\end{gather}
The solution of this problem is twice the free space Green's function of the Laplacian and given by
\begin{equation}\label{SurfaceGreens_b}
G_p(\bx;\bx_0) = \frac{1}{2\pi|\bx-\bx_0|}.
\end{equation}
\esub
If we examine the solution of \eqref{CapExtU} in the vicinity of $\bx_j$ through the stretched coordinates $\by = \eps^{-1}(\bx-\bx_j)$, we see that $G_p(\bx;\bx_j) = (2\pi \eps|\by|)^{-1}$. For this reason, we are motivated to expand the solution of \eqref{CapExtU} as 
\begin{equation}\label{OuterExp}
u = \frac{u_0}{\eps} + u_1 + \eps u_2 + \eps^2 u_3 + \bigoh(\eps^3).
\end{equation}
This implies that each of the problems for $u_j$ for $j=0,1,2,\ldots$ satisfy
\begin{equation}\label{OuterEqns}
\Delta u_j = 0, \quad \bx \in \Omega, \qquad \partial_n u_j = 0, \quad \bx \in\partial\Omega\setminus\{ \bx_1,\ldots,\bx_N\}.
\end{equation}
The solutions of \eqref{OuterEqns} are either uniform constants, or superpositions of surface Green's functions \eqref{SurfaceGreens}. The local conditions on $u_j$ as $\bx\to \bx_j$ are determined from a local solution $w(\by)$ in terms of the variable $\by = \eps^{-1}(\bx-\bx_j)$ where $\by = (s_1,s_2,\eta)$ is a local coordinate system in the vicinity of the $j^\text{th}$ pore. This local problem is expanded in a form similar to \eqref{OuterExp} 
\begin{equation}
w = \frac{w_{0,j}}{\eps} + w_{1,j} + \eps w_{2,j} + \eps^2 w_{3,j} + \bigoh(\eps^3),
\end{equation}
where in the vicinity of the $j^\text{th}$ pore, each sub-problem $w_{k,j}$ satisfies a single pore equation
\bsub\label{Innertype}
\begin{gather}
(\partial_{s_1s_1} + \partial_{s_2s_2} + \partial_{\eta\eta}) w_{k,j} = 0\,, \qquad \eta>0\,, \quad (s_1,s_2)\in\mathbb{R}^2;
\\[5pt]
w_{k,j} = 0\,, \qquad \eta = 0\,, \quad s_1^2 + s_2^2 < a_j^2\,;
\qquad \partial_{\eta} w_{k,j} = 0\,, \quad \eta = 0\,, \quad s_1^2 + s_2^2 
\geq a_j^2 \,.
\end{gather}
\esub
Each of these problems has an exact solution of form
\begin{equation}
   w_{k,j} = w_{k,j}(\infty) \left(1-w_c\right) \,,\label{Leadingv0:sol}
\end{equation}
where $w_{k,j}(\infty)$ is the constant far field solution and $w_c$ is the solution to the electrified disk problem 
\bsub\label{wc}
\begin{gather}
(\partial_{s_1s_1} + \partial_{s_2s_2} + \partial_{\eta\eta}) w_c = 0\,, \qquad \eta>0\,, \quad (s_1,s_2)\in\mathbb{R}^2;
\\[5pt]
w_c = 1\,, \qquad \eta = 0\,, \quad s_1^2 + s_2^2 < a_j^2\,;
\qquad \partial_{\eta} w_c = 0\,, \quad \eta = 0\,, \quad s_1^2 + s_2^2 
\geq a_j^2 \,.
\end{gather}
The exact solution to problem \eqref{wc} is (cf. page 38 of \cite{FA})
\begin{equation}\label{wc:sol}
w_c = \frac{2}{\pi}\sin^{-1}\Big(\frac{a_j}{L}\Big), \qquad L \equiv
 \frac12  \Big[ \sqrt{ [ (s_1^2+ s_2^2)^{\frac12} + a_j]^{2} + \eta^2} 
+ \sqrt{ [ (s_1^2+ s_2^2)^{\frac12} - a_j]^{2} + \eta^2} \Big]\,. 
\end{equation}
In terms of the capacitance $c_j={2a_j/\pi}$ of the $j^\text{th}$ pore,
 the far-field behavior
\begin{equation}\label{wc_farfield}
w_c \sim c_k \left( \frac{1}{\rho} + \frac{\pi^2 c_k^2}{24} \left(
\frac{1}{\rho^3} - \frac{3\eta^2}{\rho^5} \right) + \cdots \right)\,,
\qquad \mbox{as} \quad \rho \equiv \sqrt{s_1^2+s_2^2+\eta^2}\to\infty \,; \qquad c_k \equiv 
\frac{2a_k}{\pi}\,.
\end{equation}
\esub
is obtained. Since $\rho \sim |\bx-\bx_j|$, the far field behavior \eqref{wc_farfield} together with \eqref{Leadingv0:sol} implies that as $\bx\to\bx_j$, the matching condition with \eqref{OuterExp} yields
\[
\frac{w_{0,j}}{\eps} \sim \frac{w_{0,j}(\infty)}{\eps}\left(1- \frac{\eps c_j}{|\bx-\bx_j|} \right) \sim \frac{w_{0,j}(\infty)}{\eps} - \frac{w_{0,j}(\infty) c_j}{|\bx-\bx_j|}  \sim \frac{u_0}{\eps} + u_1 + \cdots
\]
This condition implies that $u_0 = w_{0,j}(\infty) $ for all $j=1,\ldots,N$ so that $u_0$ is a constant. It also provides a local singularity condition on $u_1$ so that it solves the problem
\bsub\label{OuterU1}
\begin{gather}
\label{OuterU1_a} \Delta u_1 = 0, \quad \bx \in \Omega; \qquad \partial_n u_1 = 0, \quad \bx \in\partial\Omega\setminus\{ \bx_1,\ldots,\bx_N\},\\[5pt]
\label{OuterU1_b} u_1(\bx) \sim \frac{-u_0 c_j}{|\bx - \bx_j|} + \cdots, \qquad \bx\to\bx_j; \qquad j = 1,\ldots,N.
\end{gather}
\esub
In terms of the Green's function $G_p(\bx;\bx_0)$ satisfying \eqref{SurfaceGreens}, the general solution of is \eqref{OuterU1} 
\begin{equation}\label{eqnU1}
u_1 = -2\pi u_0 \sum_{j=1}^N c_j G_p(\bx;\bx_j) + \chi_1,
\end{equation}
where $\chi_1$ is a constant to be found. In the summation component of \eqref{eqnU1}, each term contributes a monopole to the far field through the limiting behavior $G_p(\bx;\bx_j) \sim(2\pi |\bx|)^{-1}$ as $|\bx|\to\infty$. The normalization condition \eqref{CapExtU_c} specifies that the combined contribution should be unity, therefore
\begin{equation}\label{eqnCbar}
u_0 = \frac{-1}{N\bar{c}}, \qquad \bar{c} = \frac{1}{N}\sum_{j=1}^{N} c_j.
\end{equation}
We now proceed to the next order in $\eps$ to calculate the constant term $\chi_1$. The first step is to find the far field constant $w_{1,j}(\infty)$ of \eqref{Leadingv0:sol} near the $j^\text{th}$ hole. This comes from the local behavior of \eqref{eqnU1} as $\bx\to\bx_j$,
\bsub\label{u1asy}
\begin{equation}\label{u1asy_a}
u_1(\bx) \sim \frac{-u_0 c_j}{|\bx-\bx_j|} + B_j + \chi_1, \qquad \bx \to \bx_j, \qquad j = 1,\ldots,N,
\end{equation}
where the constants $B_j$ are given by
\begin{equation}\label{u1asy_b}
B_j = -2\pi u_0 \sum_{\substack{k=1 \\ k\neq j}}^N c_k G(\bx_j;\bx_k) = -u_0\sum_{\substack{k=1 \\ k\neq j}}^N\frac{c_k}{|\bx_j - \bx_k|}.
\end{equation}
\esub
This yields that $w_{1,j}(\infty) = B_j + \chi_1$ so that $w_{1,j} = (B_j + \chi_1)(1-w_c)$ and provides a local singularity condition on $u_2$ so that this problem now satisfies
\bsub\label{OuterU2}
\begin{gather}
\label{OuterU2_a} \Delta u_2 = 0, \quad \bx \in \Omega, \qquad \partial_n u_2 = 0, \quad \bx \in\partial\Omega\setminus\{ \bx_1,\ldots,\bx_N\},\\[5pt]
\label{OuterU2_b} u_2(\bx) \sim \frac{- c_j(B_j + \chi_1)}{|\bx - \bx_j|}, \qquad \bx\to\bx_j; \qquad j = 1,\ldots,N.
\end{gather}
\esub
In terms of the Green's function $G_p(\bx;\bx_0)$ satisfying \eqref{SurfaceGreens}, the general solution of is \eqref{OuterU2} is
\begin{equation}\label{eqnU2}
u_2(\bx) = -2\pi \sum_{j=1}^N c_j (B_j + \chi_1) G_p(\bx;\bx_j) + \chi_2,
\end{equation}
where $\chi_2$ is a constant. The normalization condition \eqref{CapExtU_c} was satisfied exactly by the equation for $u_1$, therefore \eqref{eqnU2} must make no contribution to the monopole as $|\bx|\to\infty$. This requires the solvability condition $\sum_{j=1}^N c_j (B_j + \chi_1)=0$ which fixes the value of $\chi_1$ to be
\begin{equation}\label{eqnChi1}
\chi_1 = u_0 \sum_{j=1}^{N} c_j B_j = -2\pi u_0^2 \sum_{j=1}^N \sum_{\substack{k=1 \\ k\neq j}}^N c_j c_k G_p(\bx_j;\bx_k) = - u_0^2 \sum_{j=1}^N \sum_{\substack{k=1 \\ k\neq j}}^N \frac{c_j c_k}{|\bx_j - \bx_k|}.
\end{equation}
One more application of this process is relatively simple and yields the next correction term $\chi_2$. The local behavior of \eqref{eqnU2} as $\bx\to\bx_j$ is given by
\bsub\label{u2asy}
\begin{equation}\label{u2asy_a}
u_2(\bx) \sim \frac{- c_j(B_j + \chi_1)}{|\bx-\bx_j|} + D_j + \chi_2, \qquad \bx \to \bx_j, \qquad j = 1,\ldots,N,
\end{equation}
with the constants $D_j$ given by
\begin{equation}\label{u2asy_b}
D_j = -2\pi \sum_{\substack{k=1 \\ k\neq j}}^N c_k (B_k + \chi_1)G_p(\bx_j;\bx_k).
\end{equation}
\esub
This yields that $w_{2,j}(\infty) = D_j + \chi_2$ so that $w_{2,j} = (D_j + \chi_2)(1-w_c)$ and so $u_3$ satisfies
\bsub\label{OuterU3}
\begin{gather}
\label{OuterU3_a} \Delta u_3 = 0, \quad \bx \in \Omega; \qquad \partial_n u_3 = 0, \quad \bx \in\partial\Omega\setminus\{ \bx_1,\ldots,\bx_N\};\\[5pt]
\label{OuterU3_b} u_3(\bx) \sim \frac{- c_j(D_j+\chi_2)}{|\bx - \bx_j|} + \cdots, \qquad \bx\to\bx_j; \qquad j = 1,\ldots,N.
\end{gather}
\esub
Again, $u_3$ does not contribute a monopole as $|\bx|\to\infty$, so the condition $\sum_{j=1}^Nc_j(D_j+\chi_1)=0$ must be imposed which yields that
\begin{align}
\nonumber \chi_2 &= u_0\sum_{j=1}^N c_j D_j = -2\pi u_0 \sum_{j=1}^N \sum_{\substack{k=1 \\ k\neq j}}^N c_k c_j (\chi_1+B_j)G_p(\bx_k;\bx_j)\\
{} & = \frac{\chi_1^2}{u_0} -2\pi u_0\sum_{j=1}^N \sum_{\substack{k=1 \\ k\neq j}}^N c_k c_j B_j G_p(\bx_k;\bx_j) \label{eqnChi2}
\end{align}
At this point, we recall that the goal is to determine the constant term $\lim_{|\bx|\to\infty}u(\bx) = -C_p^{-1}$ in the far field expansion \eqref{CapExtU_c}. From the expansion \eqref{OuterExp}, we have that
\[
\frac{-1}{C_p} = \frac{u_0}{\eps} \left[1 +\frac{\eps \chi_1}{u_0} + \frac{\eps^2 \chi_2}{u_0}  + \bigoh(\eps^3) \right].
\]
Rearranging this expression for the capacitance $C_p$ and simplifying yields
\begin{equation}\label{eqnC}
C_p = -\frac{\eps}{u_0}\left[1 - \frac{\eps\chi_1}{u_0} + \frac{\eps^2 \chi_1^2}{u_0^2} - \frac{\eps^2\chi_2}{u_0} + \bigoh(\eps^3) \right].
\end{equation}
Now using the expression for $u_0$ in \eqref{eqnCbar}, $\chi_1$ in \eqref{eqnChi1} and $\chi_2$ in \eqref{eqnChi2} yields that
\[
C_p = \eps N\bar{c} - \eps^2\sum_{j=1}^N \sum_{\substack{k=1 \\ k\neq j}}^N \frac{c_j c_k}{|\bx_j - \bx_k|} + \eps^3 \sum_{j=1}^N \sum_{\substack{k=1 \\ k\neq j}}^N\sum_{\substack{i=1 \\ i\neq j}}^N \frac{c_i c_j c_k}{|\bx_j - \bx_k||\bx_j - \bx_i|} + \bigoh(\eps^4).
\]
We now recall that $c_k = 2a_k/\pi$ and obtain the final simplified expression for $C_p$
\begin{equation}\label{finalC}
C_p = \frac{2N\eps\bar{a}}{\pi} - \frac{4\eps^2}{\pi^2} \sum_{k\neq j} \frac{a_j a_k}{|\bx_j - \bx_k|} + \frac{8\eps^3}{\pi^3} \sum_{\substack{j\neq k \\ i\neq j}} \frac{a_i a_j a_k}{|\bx_j - \bx_k||\bx_j - \bx_i|}+ \bigoh(\eps^4), \qquad \bar{a} = \frac{1}{N} \sum_{j=1}^{N} a_j.
\end{equation}
The corresponding flux $J_p = 2\pi DC_p$, determined from \eqref{FluxFormulaPlane}, is then given by
\begin{equation}\label{AsyPlaneFlux}
J_p = 4DN\eps\bar{a}\Big[ 1 - \frac{2\eps}{N\pi \bar{a}} \sum_{k\neq j} \frac{a_j a_k}{|\bx_j - \bx_k|} + \frac{4\eps^2}{N\pi^2\bar{a}} \sum_{\substack{j\neq k \\ i\neq j}} \frac{a_i a_j a_k}{|\bx_j - \bx_k||\bx_j - \bx_i|} + \bigoh(\eps^3) \Big].
\end{equation}
The result \eqref{mainResult} follows from setting $a_j = 1$ in \eqref{AsyPlaneFlux} for $j=1,\ldots N$.

As a check on the validity of \eqref{AsyPlaneFlux}, we compare to an exact solution in the $N=2$ case for two unit discs separated by $d = |\bx_1 -\bx_2|$.  A separable solution of \eqref{CapExt} in bi-polar coordinates \cite{Strieder08,Strieder12} determined that
\begin{equation}\label{ExactStrieder}
J_p = 8D\left( 1 -\frac{2}{\pi d} + \frac{4}{\pi^2 d^2} - \frac{2(12 + \pi^2)}{3\pi^3 d^3} + \frac{16(3+ \pi^2)}{3\pi^4 d^4} - \frac{4(120 + 70\pi^2 + 3\pi^4)}{15\pi^5 d^5} \right) + \bigoh(d^{-6}), \quad \mbox{as} \quad d\to\infty.
\end{equation}
We remark that \eqref{ExactStrieder} is a corrected version of equation (28) in \cite{Strieder08} which amends a small algebraic error carried from their previous equation (27).

By setting $\eps a_j = 1$ in \eqref{AsyPlaneFlux} for $j = 1,\ldots, N$, the first three terms of \eqref{ExactStrieder} and \eqref{AsyPlaneFlux} are in agreement.  The leading order term of \eqref{AsyPlaneFlux} is the classic Berg-Purcell result \eqref{resultBP} while the higher order terms give corrections due to pair and triplet pore interactions, respectively. This analysis can be extended to obtain further corrections to the expansion for the flux as part of a multipole expansion of solutions to \eqref{CapExtU}.

In the following section, we develop a numerical method which enables precise validation of the flux expressions \eqref{AsyPlaneFlux} and \eqref{ExactStrieder}.
 
\section{A Boundary Spectral Method for the capture problem}\label{sec:numerics}

In this section we outline a numerical spectral boundary element method for the exterior mixed Neumann-Dirichlet boundary value problems \eqref{CapExt}.  Our formulation will highlight the similarities between the numerical solution to the the planar and the spherical problem.

In our numerical method, it is convenient to solve a problem which decays as $|\bx|\to\infty$ and so we consider the equivalent capture problem for $u(\bx) = 1 - v(\bx)$ where $v(\bx)$ satisfies \eqref{CapExt} and $u(\bx)$ solves
\begin{equation}\label{MainRescaled_a}
 \Delta u = 0\,, \quad \bx\in\Omega\,;
 \qquad u = 1\,, \quad \bx \in \Gamma_a\,, \qquad \partial_n u = 0\,, \quad \bx \in\Gamma_r,
\end{equation}
where $\Gamma_a$ is a set of absorbing circular pores and $\Gamma_r$ is the reflecting complement of $\Gamma_a$ on the boundary $\partial\Omega$. We complete the problem by specifying that the solution of $u(\bx)$ decays to zero in the far field. In this formulation, the capacitance $C$ of the target is specified by the flux $J$ over $\Gamma_a$ and satisfies
\begin{equation}
 J \equiv \int_{\Gamma_a} \partial_n u\, \dS =
 \begin{cases}
 2 \pi C & \text{Planar Capture}\\
 4 \pi C & \text{Spherical Capture}
 \end{cases}
\end{equation}
which, together with the decay condition, determines the far-field behavior
 \begin{equation}\label{MainRescaled_b}
u(\bx) = \frac{C}{|\bx|} + \mathcal{O}\left(\frac{1}{|\bx|^2}\right) , \qquad |\bx|\to\infty .
\end{equation}

We formulate the numerical problem as a linear integral equation, specifically a \emph{Neumann to Dirichlet map} \cite{Greenbaum93,Tausch1998} on the set of pores, $\Gamma_a$, relating the known surface potential, $u|_{\partial\Omega} = p(\bx)$, equal to unity in $\Gamma_a$  to 
the surface flux, $\partial_n u|_{\pOm}  = q(\bx)$, 
which is unknown on $\Gamma_a$ and vanishes on $\Gamma_r$. Fortunately, the exact solution to the Neumann problem is known in terms of the surface Green's functions \cite{K,FA,sneddon1966mixed,SurfaceGreen3D}
$$u(\bx) = \int_{\by \in \pOm}  G(\bx;\by) \, q(\by)~ \dS, \qquad \bx\in \Om .$$
We simplify this by first noting that the surface flux, $q(\bx)$, is non-zero only on the pores, $\Gamma_a$, and second by restricting our interest to the surface  where $u(\bx)=p(\bx)$. This yields the linear integral equation
\begin{equation}\label{PoreIE}
p(\bx)  = \A \left [ q(\bx) \right] \equiv  \frac{1}{2\pi} \int_{\by \in \Gamma_a}  g\left (\left | \bx-\by \right | \right ) q(\by)~ \dS, \qquad \bx \in
\pOm , 
\end{equation}
where the kernel of the integral operator is defined by the Green's function restricted to the surface
$$G(\bx;\by) = \frac{1}{2\pi} g\left (\left | \bx-\by \right | \right ) \qquad \text{for} \ \bx,\by \in \pOm $$
where
\begin{align}
 g(\mu) &= g_p(\mu) \equiv \frac{1}{\mu} \phantom{+ \frac{1}{2}  \log\left(\frac{\mu}{2+\mu}\right) }\qquad \text{(Planar Capture)}\\
 g(\mu) &= g_s (\mu)\equiv \frac{1}{\mu} + \frac{1}{2}  \log\left(\frac{\mu}{2+\mu}\right) \qquad \text{(Spherical Capture)}
 \end{align}
as defined in \eqref{SurfaceGreens} and \eqref{FluxSphere_b} respectively.  We are interested in the specific case when the surface potential $p(\bx)=1$ for $\bx \in \Gamma_a$.

We solve this problem pseudo-spectrally by a judicious choice of  basis of functions for the surface potential, $p(\bx)$, and the surface flux, $q(\bx)$, within the pores $\Gamma_a$.   We are guided by the known exact solution \cite{FA, sneddon1966mixed} for a single absorbing circular pore  on a half plane \eqref{wc}.   At this point we will simplify the calculation by assuming that the $\Np$ pores, $\Gamj$,  that constitute $\Gama = \displaystyle\cup_{j=1}^\Np \Gamj$ have a common radius $\alpha$.  For the planar problem the $\Np$ pores are discs of radius $\alpha$ centered at  points $\left \{ \bx_1, \ldots, \bx_\Np \right \}$ on the plane $z=0$. 

For the spherical problem the $\Np$ pores that each subtend an angle $\sig$ centered at  points $\left \{ \bx_1, \ldots, \bx_\Np \right \}$. The boundary of each pore is the set of points on the sphere which are a distance $\alpha= 2\sin (\sig/2)$ from its center. On the surface near each pore, we introduce a local spherical coordinate system $(\theta_k, \phi_k)$ with the polar axis aligned with the pore's center. We observe that if we make a change of variables
$$\xi_k = 2 \sin \theta_k/2 , \qquad  t_k=\phi_k, $$
that the $k^\text{th}$ pore occupies a disc in $(\xi_k, t_k)$ space,
$$\Omega_k = \left \{( \xi_k,t_k) \ | \ 0 \le \xi_k \le \alpha, \ 0 < t_k \le 2 \pi \right \}  \qquad \text{where} \quad \alpha= 2\sin (\sig/2),$$
and the spherical area element can be rewritten as
$$dS = \sin (\theta_k) \ d \theta_k \, d\phi_k = \xi_k d \xi_k \, d t_k, $$
which is identical to the area element for planar polar coordinates. 

\begin{figure}
\centering
\includegraphics[width = \textwidth]{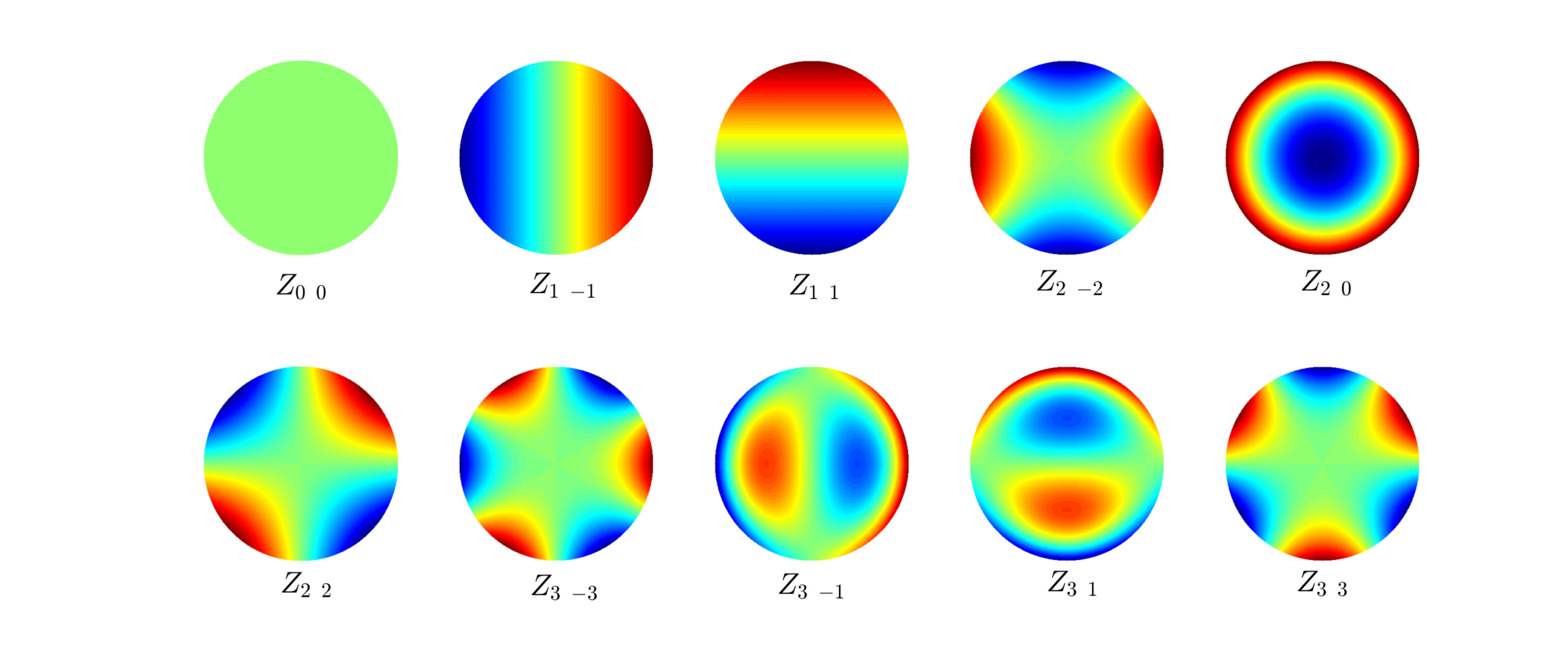}
\parbox{0.95\textwidth}{\caption{ The first $10$ Zernike polynomials on the unit disc.  \label{fig:First10Zern}}}
\end{figure}

We will now choose a basis for the surface potential $p(\bx)$ on each pore.
If we define a Cartesian coordinate system $\left (X,Y \right )=\left (\xi  \cos (t), \xi \sin(t) \right )$ for a given pore, a natural basis would be polynomials in $(X,Y)$ of degree less than or equal to $M$. An orthonormal basis for these polynomials on the unit disc (expressed in polar coordinates) are \emph{Zernike Polynomials} \cite{Zernike2013}, defined by
$$Z_{m\,j} (\xi, t) = 
\begin{cases} 
 P_{m\,j}(\xi) \sin(jt) & j>0 , \\
 P_{m\,0}(\xi) & j=0, \\
 P_{m\,|j|}(\xi) \cos(jt) & j<0,
\end{cases}
 \qquad m=0,1, \ldots ,M, \quad j=-m, -m+2, \ldots ,m-2, m .
$$
Each $P_{m\,j}(\xi)$ is a degree $m$ polynomial containing terms of degree $j,j+2, \cdots,m-2,m$. The first $Z_{m\, j}$ are
\begin{gather}
\nonumber Z_{0\, 0} = \frac{1}{\sqrt{\pi}} ,
\\
\nonumber Z_{1\, -1}= \frac{2}{\sqrt{\pi}}  \xi \cos  t , \quad 
Z_{1\, 1}= \frac{2}{\sqrt{\pi}} \xi \sin t,
\\
\nonumber Z_{2\, -2}= \sqrt{\frac {6}{\pi}}\xi^2 \cos 2t , \quad
Z_{2\, 0}= \sqrt{\frac {3}{\pi}}\left (2\xi^2-1\right ) , \quad
Z_{2\, 2}= \sqrt{\frac {6}{\pi}}\xi^2 \sin 2t , \\
\nonumber Z_{3\, -3}=\sqrt{\frac {8}{\pi}}\xi^3 \cos 3t  , \quad
Z_{3\, -1}=\sqrt{\frac {8}{\pi}}\left (3\xi^3-2\xi \right ) \cos t , \quad
Z_{3\, 1}=\sqrt{\frac {8}{\pi}}\left (3\xi^3-2\xi \right ) \sin t , \quad
Z_{3\, 3}=\sqrt{\frac {8}{\pi}}\xi^3 \sin 3t,
\end{gather}
and plotted in Fig.~\ref{fig:First10Zern}.
There are $m+1$ polynomials of degree $m$ and a total of $(M+1)(M+2)/2$ polynomials of degree $M$ or less. If we define the inner-product on the disc of radius $\alpha$, denoted by $\Omega$ here,
$$\left \langle \Phi(\xi,t) , \Psi(\xi,t)  \right \rangle_\Omega
\equiv \int_{t=0}^{2\pi}\int_{\xi=0}^\alpha \Phi(\xi,t)  \Psi(\xi,t)\,\xi \,d \xi \, d t,
$$
the orthonormality condition for the Zernike polynomials on the unit disc ($\alpha=1$) is 
$$\left \langle Z_{m\, j}(\xi,t), Z_{m'\,j'}(\xi,t) \right \rangle_\Omega = \delta_{m\,m'} \delta_{j\,j'} .$$
For discs of radius $\alpha$, one can use a rescaled basis of Zernike polynomials, $ Z_{m\,j}(\xi/\alpha,t)$
for which the orthogonality condition reads
$$\left \langle Z_{m\,j}(\xi/\alpha,t), Z_{m'\,j'}(\xi/\alpha,t) \right \rangle_\Omega = \alpha^2 \delta_{m\,m'} \delta_{j\,j'} .$$

We can now approximate the known surface potential on the $k^\textrm{th}$ disc, $\Gamma_k$, as a linear combination of the Zernike polynomials up to degree $M$
\begin{equation}
\label{PotNumExp}
p(\xi_k,t_k) = \sum_{m=0}^M \sum_{\ell=0}^{m}  {c_{m\,j\,k}}Z_{m\,j}(\xi_k/\alpha,t_k),  \qquad j =2\ell -m,
\end{equation}
where the total number of coefficients is $\Np (M+1)(M+2)/2$, the number of pores multiplied by the number of polynomials of degree $M$ or smaller.
If $p(\xi_k,t_k) =1$, we find that for each pore there is a single non-zero mode with $m=j=0$, 
\begin{equation}
\label{PotNumVal}
c_{m\,j\,k} = \alpha^2 \sqrt{\pi} \, \delta_{m\,0}\, \delta_{j\,0},
\end{equation}
for each pore $\Omega_k$.

We now need to find a basis to approximate the flux on each of the pores. It is well known that the flux may be singular for mixed Neumann-Dirichlet problems \cite{sneddon1966mixed,Duffy08} and this problem is not an exception. We appeal to an analogous problem, specifically a single circular pore on a plane bounding a half-space, for insight. This problem can also be formulated as an integral equation of the form \eqref{PoreIE}. The kernel in the planar case is equivalent to the most singular term in the spherical problem, $g(\mu) =1/\mu$, and an exact solution is given in \eqref{wc}. If the surface potential is of the form of a Zernike polynomial 
$$p(\xi,t) =  Z_{m\,j} (\xi/\alpha, t) = 
\begin{cases} 
 P_{m\, j}(\xi/\alpha) \sin(jt) & j>0 ; \\
 P_{m\,0}(\xi/\alpha) & j=0; \\
 P_{m\,|j|}(\xi/\alpha) \cos(jt) & j<0,
\end{cases}
$$
an exact solution for the flux in the pore can be found of the form 
$$q(\xi,t) =\frac{1} {\sqrt{\alpha^2-\xi^2}}\begin{cases} 
 Q_{m\,j}(\xi/\alpha) \sin(jt) & j>0 ; \\
 Q_{m\,0}(\xi/\alpha) & j=0; \\
 Q_{m\,|j|}(\xi/\alpha) \cos(jt) & j<0,
\end{cases} 
$$
where $Q_{m\,j}$ is a polynomial of degree $m$ containing terms of degree $j,j+2, \cdots,m-2,m$. In the reflecting region exterior to the pore ($\xi > \alpha$), the flux vanishes identically.  An important example of this exact solution is the constant surface potential in a single pore where for $p(\xi,t)=\sqrt{\pi}Z_{0\, 0} =1$, the surface flux is 
$$q(\xi,t)=  \frac{2} {\pi \sqrt{\alpha^2-\xi^2}},\qquad 0\leq\xi<\alpha,$$
which is positive and exhibits an inverse square root singularity at the edge of the pore  (which is also evident in \eqref{wc}). The singularity is integrable and the flux can be computed as
$$ J = \int_{\xi=0}^\alpha \int_{t=0}^{2 \pi} q(\xi,t) \, \xi\, d\xi \, dt = {4 \alpha}, $$
which allows us to recover the backbone of the Berg-Purcell result \eqref{resultBP}, namely that the flux for a single pore (appropriately non-dimensionalized) is four times the perimeter.

This exact solution suggests that an appropriate basis for the surface flux is terms of the form
$$q_{m\,j}(\xi,t) \equiv  \frac{1} {\sqrt{\alpha^2-\xi^2}}\begin{cases} 
 (\xi/\alpha)^m \sin(jt) & j>0;  \\
 (\xi/\alpha)^m & j=0; \\
 (\xi/\alpha)^m \cos(jt) & j<0,
\end{cases}
\qquad m=0,1, \ldots ,M, \quad j=-m, -m+2, \ldots ,m-2, m .
$$
This basis captures the nature of the flux singularity on the boundary of the pore and spans the exact solution for the single pore problem on the half plane with a polynomial flux function of degree $M$ or less. Our assumption (borne out by the asymptotics of \S \ref{sec:planarpore} and in \cite{LWB2017} and our numerics) is that the corrections to the flux due to the curvature of the surface and the pore interactions are subdominant and can also be captured by this basis.

As such we expand the surface flux on the $k^{\text{th}}$ pore as a sum of these functions,
\begin{equation}
\label{FluxNumExp}
q(\xi_k,t_k) = \sum_{m=0}^M \sum_{\ell=0}^m b_{m\,j\,k} q_{m\,j}(\xi_k,t_k), \qquad j=2\ell-m, 
\end{equation}
where the $N(M+1)(M+2)/2$ constants $b_{m\,j\,k}$ need to be determined. Substituting the expansion for the surface potential \eqref{PotNumExp} and the surface flux \eqref{FluxNumExp}  into the governing integral equation \eqref{PoreIE} yields
$$ \sum_{k'=1}^\Np\sum_{m'=0}^M \sum_{\ell'=0}^{m'}   {c_{m'j'k'}}Z_{m'j'}(\xi_{k'}/\alpha,t_{k'})= \sum_{k=1}^{\Np} \sum_{m=0}^M \sum_{\ell=0}^m b_{m\,j\,k}  \A[q_{m\,j}(\xi_k,t_k)] +\E_M, \qquad  j'=2\ell'-m', \quad j=2\ell-m,$$
where $\E_M$ is the error incurred by having a finite approximation of order $M$ for both the surface flux and the surface potential. We now project both sides of the equation onto the Zernike polynomial basis for the flux functions; applying 
the operator $ \left \langle Z_{m'j'}(\xi_{k'}/\alpha,t_{k'}), \cdot \right \rangle_{\Omega_{k'}}$ yields
\begin{equation}\label{BigLin}
c_{m'j'k'} = A_{m'j'k'mjk} \, b_{m\,j\,k}, \qquad  A_{m'j'k'mjk}=\frac{1}{\alpha^2} \left \langle Z_{m'j'}(\xi_{k'}/\alpha,t_{k'}), \A[q_{m\,j}(\xi_k,t_k)] \right \rangle_{\Omega_{k'}},
\end{equation}
where we have divided by $\alpha^2$. The resulting square linear system is 
of size $\Np(M+1)(M+2)/2$ for the unknown surface flux coefficients  $b_{m\,j\,k}$ in terms of the 
known surface potential coefficients $c_{m'j'k'}$, evaluated above in \eqref{PotNumVal}.
A solution to this linear system will minimize the $L^2$ norm of the error $\E_M$ over the collection of pores, $\Gamma_a$.

To evaluate the coefficients $A_{m'j'k'mjk}$ naively one needs to evaluate a quadruple integral, integrating over the discs $\Omega_k$ and $\Omega_{k'}$. However, the symmetries of the problem simplifies these evaluations immensely. First, we evaluate the surface potential induced by $q_{m\,j}(\xi,t)$,
\begin{equation}\label{eqn:surfacePotential}
p_{m\,j}(\xi,t)  \equiv  \A[q_{m\,j}(\xi,t)].
\end{equation}
This function $p_{m\,j}$ has the same angular ($t$) dependence as $q_{m\,j}$ and the $\xi$ dependence is computed numerically and tabulated for each value of $m$ and $j$ to allow for later interpolation. We now discuss some implementation details for the method, treating the planar and spherical cases separately.

\underline{Case I (Plane)}: Here $g(\mu) = 1/\mu$ and the potential \eqref{eqn:surfacePotential} induced by $q_{m\, j} = \cos(j t)\, (\xi/\alpha)^m\, (\alpha^2 - \xi^2)^{-\frac12}$ for $j\ge 0$ is (for $j<0$ replace $\cos (jt)$ by $\sin(jt)$ throughout)
\begin{align}
\nonumber p_{m\, j}(\xi,t) &= -\frac{\alpha^{-m}}{2\pi} \ds\int_{\rho=0}^{\alpha} \int_{\eta = 0}^{2\pi} \frac{\rho^m \cos (j\eta) }{ \sqrt{\alpha^2 - \rho^2} } \frac{1}{\sqrt{ \rho^2 + \xi^2 - 2 \rho \xi \cos(t-\eta) } }\, \rho\, d\rho\, d\eta\\[5pt]
\label{eqn:pjm} &= - \cos(jt) \frac{ \alpha^{-m} }{2\pi} \int_{\rho = 0}^{\alpha} \frac{\rho^m }{\sqrt{\alpha^2 - \rho^2}}  H_j(\xi/\rho)\, d\rho
\end{align}
where the function $H_j(\beta)$ is defined as
\begin{equation}\label{eq:Hj}
H_j(\beta) = \int_{\tau = 0}^{2\pi} \frac{ \cos(j \tau) }{\sqrt{\beta^2 + 1 - 2 \beta \cos \tau}}\, d \tau, \qquad \beta\geq0.
\end{equation}
The numerical evaluation of the integral $H_j(\beta)$ is simplified by noting that 
\[
H_j(\beta) = \frac{1}{\beta} H_j\left( \frac{1}{\beta}\right), \qquad \beta \neq0,
\]
which restricts computations to the range $0\leq\beta\leq1$. The integral \eqref{eq:Hj} has a logarithmic singularity at $\tau =0$ as $\beta\to1$. Effective numerical evaluation of $H_j(\beta)$ in light of this singularity is aided by writing
\begin{equation}\label{eq:Hbreakdown}
 H_j(\beta) = \int_{\tau = 0}^{2\pi} \frac{ \cos(j \tau) -1}{\sqrt{\beta^2 + 1 - 2 \beta \cos \tau}}\, d \tau + \int_{\tau = 0}^{2\pi} \frac{1}{\sqrt{\beta^2 + 1 - 2 \beta \cos \tau}}\, d \tau.
\end{equation}
The first integral in \eqref{eq:Hbreakdown} is bounded and readily approximated while the second term captures the logarithmic singularity and is expressed as an elliptic integral and evaluated with the {\tt MATLAB} function {\tt ellipke}. Returning to the integral \eqref{eqn:pjm} for the surface potential and setting $\rho = \alpha \sin s$, we have that 
\begin{equation}\label{eqn:Pquadgk}
p_{m\,j}(\xi,t) = -\frac{ \cos jt}{2\pi} \int_{s=0}^{\frac{\pi}{2}} [\sin s]^m\,  H_j \Big( \frac{\xi}{\alpha\sin s} \Big)\, ds.
\end{equation}
If $\xi>\alpha$, then the integrand is bounded, however, for $\xi<\alpha$ the integrand has a logarithmic singularity at $s^{\ast} = \sin^{-1} (\xi/\alpha)$. In the later case, we split the integration interval at $s = s^{\ast}$ so that  the integrable singularity is placed on the endpoints. This integral is then evaluated with built in {\tt MATLAB} quadrature routines, notably {\tt quadgk} \cite{Shampine2008} which accommodate integrands with logarithmic boundary singularities.

The final step in the formation of the linear system \eqref{BigLin} requires the evaluation of the integrals
\begin{equation}\label{eq:LinSm}
A_{m'j'k'mjk}=\frac{1}{\alpha^2} \left \langle Z_{m'j'}(\xi_{k'}/\alpha,t_{k'}), p_{m\,j}(\xi_k,t_k) \right \rangle_{\Omega_{k'}},
\end{equation}
which represent the inner-products of $Z_{m'j'}(\xi_{k'}/\alpha,t_{k'})$ with $p_{m\,j}(\xi_k,t_k) $. For the case $k=k'$ the integral vanishes unless $j=j'$ and the angular portion can be evaluated exactly in this case reducing the problem to one dimension and we use built in {\tt MATLAB} quadrature routines \cite{Shampine2008}. For $k\neq k'$, we use a polar collocation grid on the disc $\Omega_{k'}$ with equally spaced and weighted points in the angular variable and a radial grid that is equally spaced in the square of the radial distance weighted by a $10$-point Newton-Cotes formula. This reduces each of the inner products \eqref{eq:LinSm} to a dot product of a weighted vector on the collocation points with the function $p_{m\,j}(\xi_k,t_k) $ evaluated via interpolation on the collocation points. This step can be easily parallelized over each of the matrix entries of \eqref{eq:LinSm}.

\underline{Case II (Sphere)}: We first determine the potential induced by the flux from \eqref{eqn:surfacePotential}. For two surface points 
$\bx = (\sin\theta\cos\phi,\sin\theta\sin\phi,\cos\phi)$ and $\bx' = (\sin\theta'\cos\phi',\sin\theta'\sin\phi',\cos\phi')$, the surface distance $d$ is 
\begin{equation}\label{eq:sphereD}
d = |\bx -\bx'| = \sqrt{2 - 2\sin\theta \sin \theta' \cos(\phi - \phi')  - 2\cos\theta \cos\theta' }.
\end{equation}
To reduce \eqref{PoreIE} from an integral over a spherical region to a circular region, the transformations
\begin{equation}
\xi = 2\sin(\theta/2), \qquad \eta = 2 \sin(\theta'/2),
\end{equation}
are applied such that $\xi\in[0,2]$, $\eta \in[0,2]$ and the surface distance \eqref{eq:sphereD} becomes 
\begin{equation}\label{eqn:surfaceD2}
d^2 = \xi^2 + \eta^2 - \frac{1}{2} \xi^2 \eta^2 - 2 \eta\, \xi \sqrt{1-\xi^2/4} \sqrt{1- \eta^2/4}\ \cos \tau,
\end{equation}
for $\tau = \phi - \phi'$. The integral \eqref{eqn:surfacePotential} can now be evaluated as 
\begin{align}
\nonumber p_{m\, j} &= \frac{\alpha^{-m}}{2\pi} \int_{\eta = 0}^{\alpha} \int_{\tau=0}^{2\pi} \frac{\eta ^m \cos j t }{\sqrt{\alpha^2 - \eta^2}} \left( \frac{1}{d} + \frac{1}{2}  \log \left[ \frac{d}{2+d} \right]\right)\, \eta\,d\eta\, d\tau \\[5pt]
\label{eqn:SpherePotential} & = \alpha \int_{s=0}^{\frac{\pi}{2}} [\sin s ]^{m+1} H_j(\alpha \sin s, \xi)\, ds, 
\end{align}
where in the final step, the substitution $\eta = \alpha \sin s$ was used. In this case the function $H_j$ is
\begin{equation}\label{eqn:HSphere}
H_j(\eta, \xi) = \frac{1}{2\pi} \int_{\tau = 0}^{2\pi} \cos j \tau \left( \frac{1}{d} + \frac{1}{2} \log d - \frac{1}{2}\log(2+d) \right) \, d\tau
\end{equation}
where $d = d(\eta,\xi,\tau)$ is given in \eqref{eqn:surfaceD2}. As in the planar case, the function $H_j(\eta,\xi)$ in \eqref{eqn:HSphere} has a singular integrand and must be treated with care to obtain an accurate numerical evaluation. In the decomposition 
\begin{equation}\label{eqn:HSphere2}
H_j(\eta, \xi) = \frac{1}{2\pi} \int_{\tau = 0}^{2\pi} \left[(\cos j \tau - 1) \Big( \frac{1}{d} + \frac{1}{2} \log d \Big) - \cos j\tau \log(2+d) \right] d\tau + \frac{1}{2\pi} \int_{\tau = 0}^{2\pi} \left(\frac{1}{d} +\frac{1}{2} \log d \right)  d\tau,
\end{equation}
the first integral has a bounded integrand and is readily evaluated while the second integral has a singular integrand. The singular component arising from the $1/d$ term is expressed in terms of an elliptical integral while the integral of the term $\frac{1}{2}\log d$ can be evaluated exactly from the identity (cf.~\cite{GR}),
\[
\int_0^{\pi} \log (a + b \cos x)\, dx = \pi \log \left[ \frac{a + \sqrt{a^2 - b^2}}{2}\right], \qquad a \geq |b| > 0.
\]
The values of $H_j(\eta,\xi)$ are tabulated over a grid of $(\eta,\xi)$ points for a range of $j$ and stored for the computation of the surface potential \eqref{eqn:SpherePotential}. For values $\xi>\alpha$, the integral \eqref{eqn:SpherePotential} is well behaved and easily evaluated. For $\xi<\alpha$ an integrable singularity is present at $s^{\ast} = \sin^{-1}(\xi/\alpha)$ which is resolved by dividing the integration interval at $s= s^{\ast}$ so the subsequent integrals have boundary singularities and evaluated in {\tt MATLAB} \cite{Shampine2008}.

The final step is to obtain the entries of the matrix $A$ in \eqref{BigLin} by calculating the projection of the surface potential onto the Zernike modes,
$$\left \langle Z_{m'j'}(\xi_{k'}/\alpha,t_{k'}), p_{m\,j}(\xi_k,t_k) \right \rangle_{\Omega_{k'}}.$$
To perform each integration, we first translate the $k'^{\text{th}}$ pore to the north pole followed by application of the collocation method discussed in the planar case.

Once the matrix is built, we solve the linear system \eqref{BigLin} using the {\tt MATLAB} built in matrix solver which yields the unknown weights of the flux functions, $b_{m\,j\,k}$. These weights allow us to compute the flux through each pore and in turn the total flux and capacitance of a given configuration.

This algorithm appears robust, although it has its limitations some of which we explore in the next section. Typically the quadratures are evaluated to obtain absolute errors of $10^{-15}$ and relative errors of $10^{-8}$ although these numbers may be degraded to $10^{-10}$ and $10^{-6}$ in the immediate neighborhood of a singularity. For the singular integrals, we move the boundary points inward by {\tt MATLAB}'s machine epsilon (roughly $10^{-15}$) to avoid overflows. Increasing the number of Zernike modes yields consistent answers with relative errors of about $10^{-8}$ which appears to be in part due to accumulated round-off errors.

We also note heuristically that there are two reasons this algorithm converges. First, for the biologically relevant case of pores whose separation is large compared to their radius, an expansion with Zernike polynomials up to degree $M$ effectively captures an $M^{\text{th}}$ order asymptotic approximation of the solution akin to the analysis of Section 2. Second, even for closely spaced pores we are minimizing the $L^2$ error for a degree $M$ polynomial approximation of the surface potential on the pores. We investigate this convergence below. 
 
\section{Numerical Results}\label{sec:results}

In this section we detail numerical results for the planar and spherical case. In practice, we have run with polynomials of degree up to $M=20$ for numbers of pores $\Np \le 20$ and run up to $N=2001$ pores for lower approximations ($M\le 6$). The calculations take from a minute to a few hours on a standard desktop computer. The method appears to be effective and accurate for small, widely separated pores which is the relevant asymptotic and biological limit. Accuracy is degraded if pore boundaries are nearly touching (which necessitates larger values of $M$ to resolve).

In the following examples, we benchmark the numerical accuracy by evaluating the relative error
\begin{equation}\label{RelTol}
\mathcal{E}_{\textrm{rel}}[J] = \left| \frac{J_{\textrm{num}} -  J_{\textrm{asy}}}{J_{\textrm{num}}} \right|.
\end{equation}
In the following examples, $J_{\textrm{asy}}$ is obtained from asymptotic, numerical and exact expressions for the flux $J$. 

\subsection{Planar Case}\label{resultsPlane}

In the following examples, we demonstrate the convergence of the numerical method as the number of Zernike modes $M$ increases and verify the accuracy against the asymptotic formula \eqref{AsyPlaneFlux}.

\subsubsection{Example: Two Planar Pores}\label{Ex:TwoPore}
\begin{figure}[h]
\centering
\subfigure[t][Two Pore Schematic]{\includegraphics[width=0.475\textwidth]{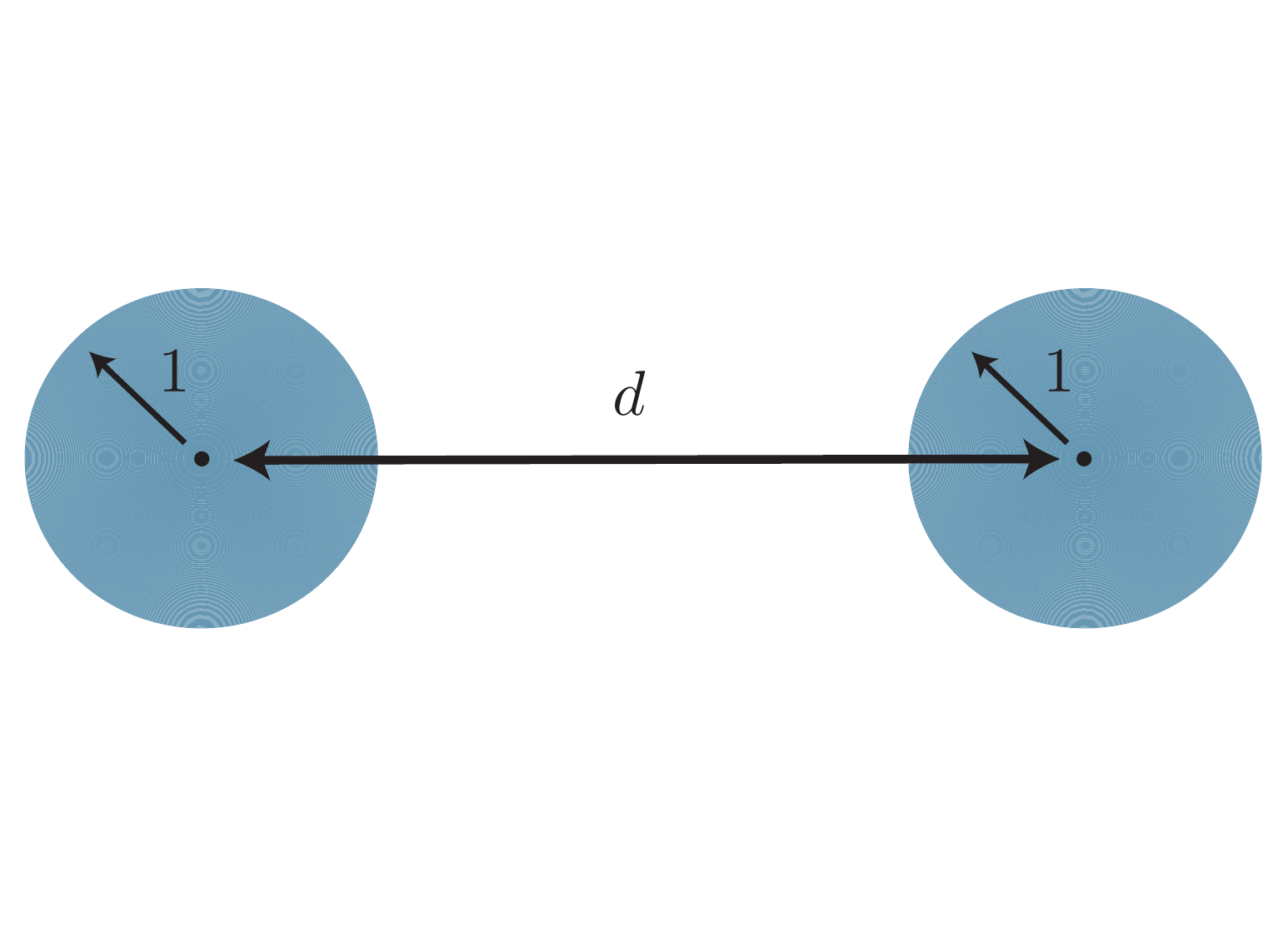} \label{fig:TwoPoreA}} \quad 
\subfigure[Convergence Results]{\includegraphics[width=0.475\textwidth]{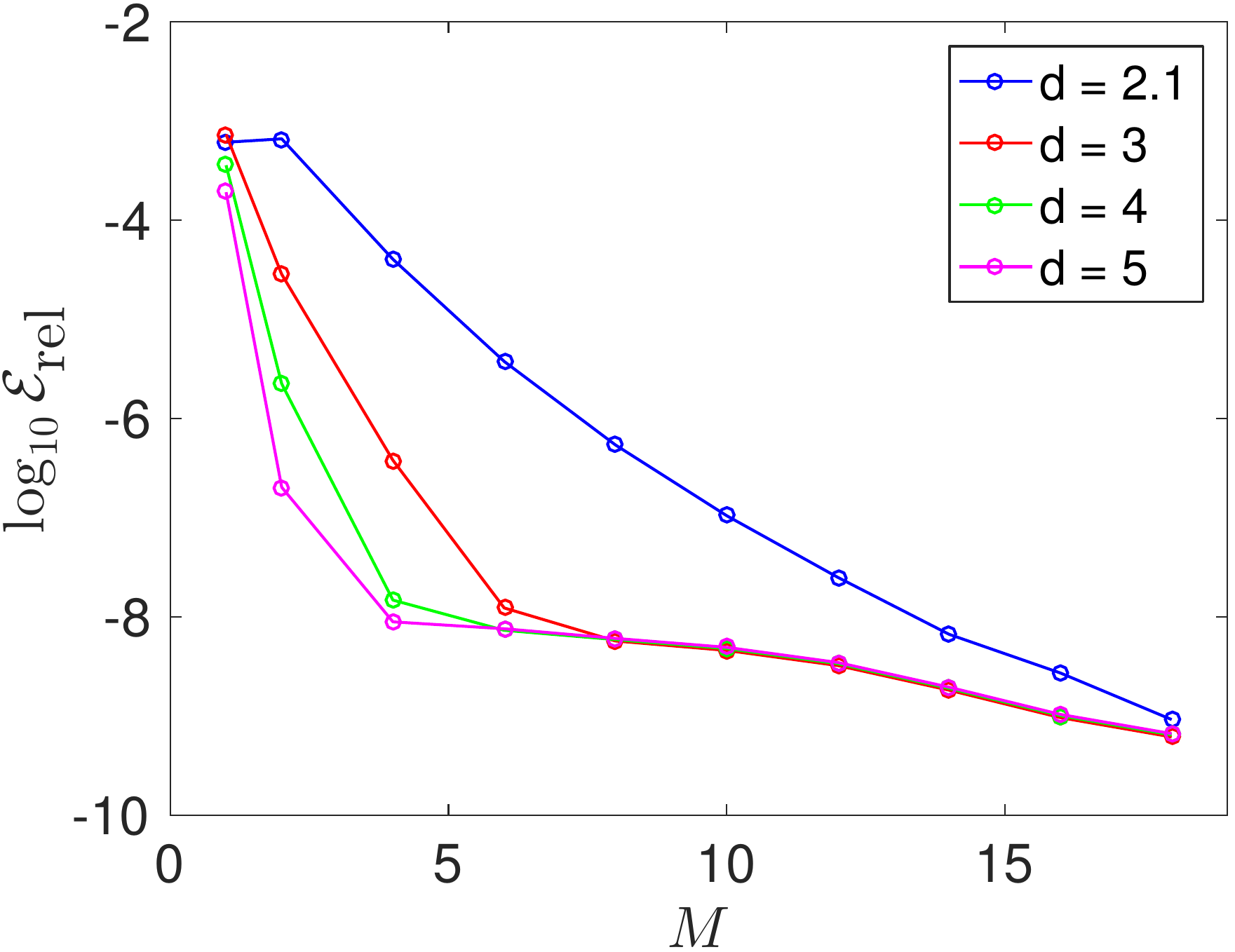} \label{fig:TwoPoreB}}
\caption{Results of example \ref{Ex:TwoPore}. \ref{fig:TwoPoreA}: Schematic of two pores with common radius $a =1$ centered at $(\pm d/2,0,0)$ and separated by distance $d$.  \ref{fig:TwoPoreB}: Convergence of the numerical relative error as the number of Zernike modes $M$ increases. When sufficient modes are included, the method has a relative error of around $10^{-8}$. Relative errors are calculated from \eqref{RelTol} with respect to a \lq\lq true\rq\rq\ solution obtained with $M=20$ modes. \label{fig:TwoPore}}
\end{figure}

In this example we take two pores of unit radius centered at $\bx = (\pm d/2,0,0)$ and demonstrate convergence of the numerical method over separation distances $d>2$. In the results of Fig.~\ref{fig:TwoPore}, we use the numerical solution for $M=20$ modes as an exact solution in \eqref{RelTol}.

The key observation from this example is that relatively few ($M\approx 6$) modes are required to accurately resolve the capture rate, provided the pore spacing is not too small. As the pore separation decreases ($d\to2^{+}$), additional modes must be included to accurately resolve the solution.

For two pores with centers separated by a distance $d$, Strieder (cf.~\cite{Strieder08,Strieder12}) calculated a series approximation for $J_p$ from a separable solution in bi-polar coordinates. The first few terms of that series and its truncation error are
\begin{equation}\label{ExactStriederResults}
J_p = 8D \left[ 1 -\frac{2}{\pi d} + \frac{4}{\pi^2 d^2} - \frac{2(12 + \pi^2)}{3\pi^3 d^3} + \frac{16(3+ \pi^2)}{3\pi^4 d^4} - \frac{4(120 + 70\pi^2 + 3\pi^4)}{15\pi^5 d^5} \right] + \bigoh(d^{-6}), \qquad d\to\infty.
\end{equation}
As remarked after equation \eqref{ExactStrieder}, the expression \eqref{ExactStriederResults} is a corrected version of equation (28) in \cite{Strieder08}.

\begin{figure}[h]
\centering
\subfigure[t][Two Pore Competition]{\includegraphics[width=0.475\textwidth]{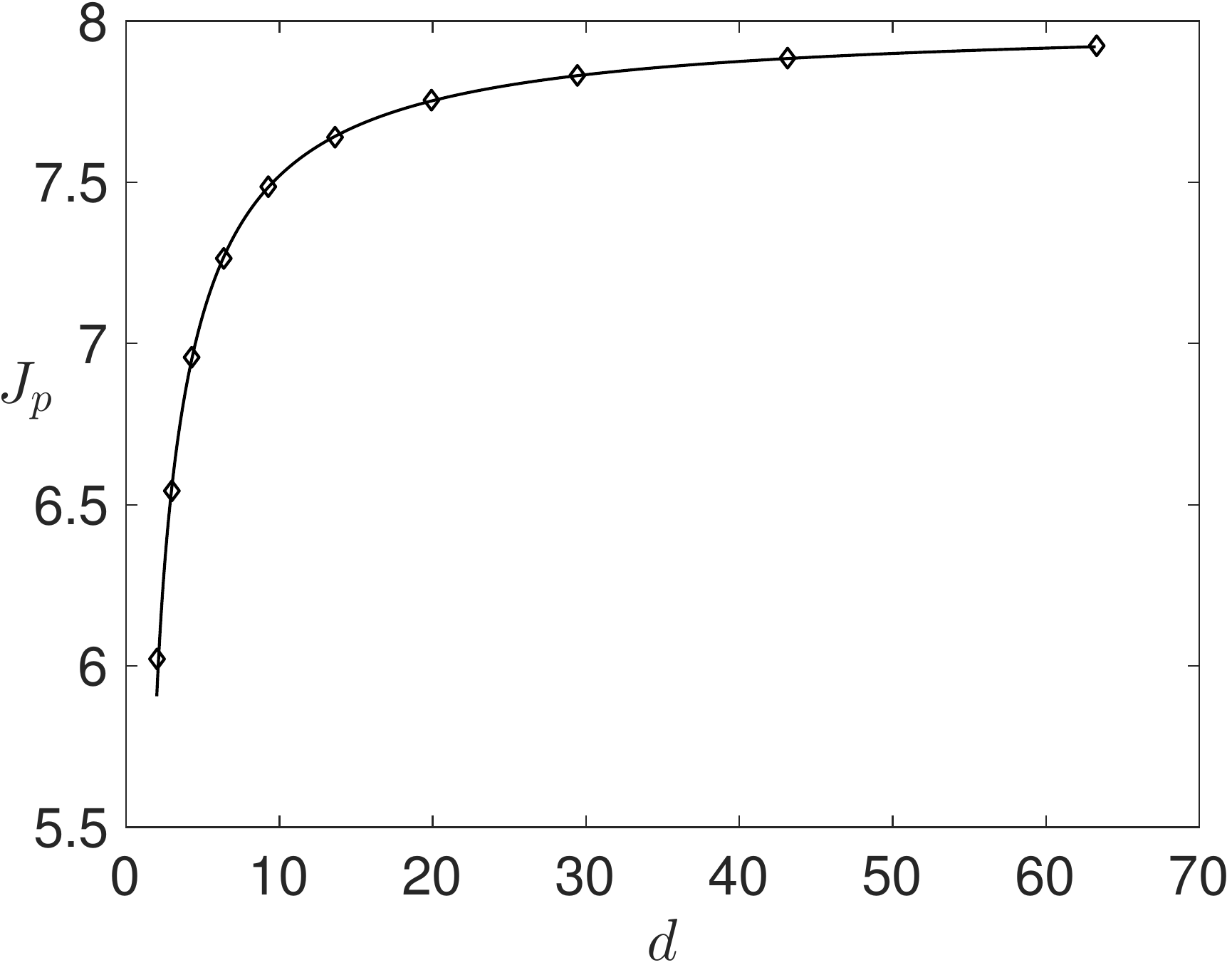} \label{fig:TwoPoreSepA}} \quad 
\subfigure[Convergence Results]{\includegraphics[width=0.475\textwidth]{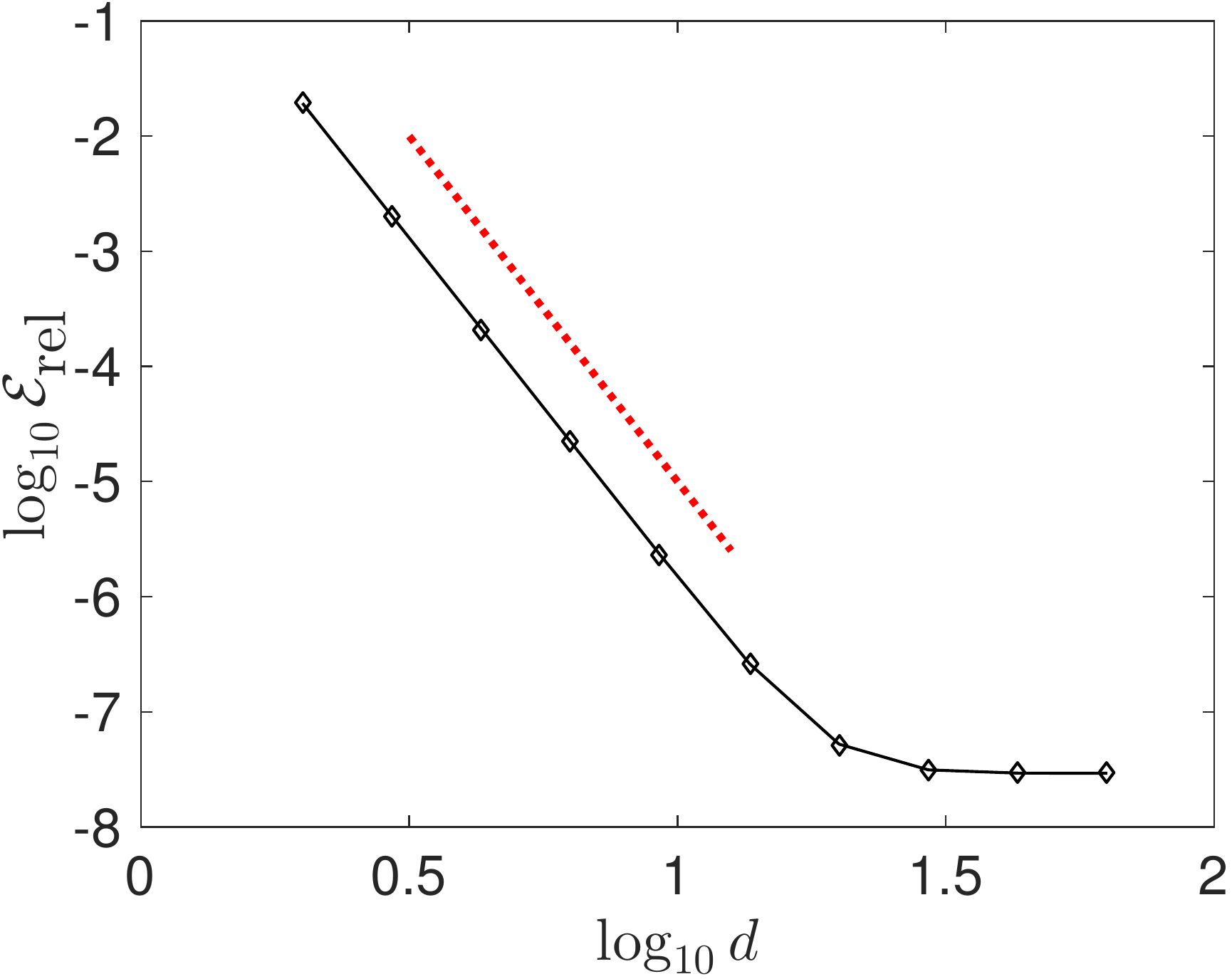} \label{fig:TwoPoreSepB}}
\caption{Results of example \ref{Ex:TwoPore} for two well separated pores. Fig.~\ref{fig:TwoPoreSepA}: The flux $J_p$ to pores with unit radius centered at $(\pm d/2,0,0)$ as given by the series (solid line) formula \eqref{ExactStriederResults} and numerical simulations (diamonds) with $M=20$ modes. At small separations $d$, interpore competition reduces the flux considerably. Fig.~\ref{fig:TwoPoreSepB}: Convergence of $\mathcal{E}_{\mathrm{rel}}$ with $M=20$ modes as the distance $d$ increases. The line (dotted red) of slope $-6$ confirms the accuracy of the series solution \eqref{ExactStriederResults}. The method accurately resolves the flux to one part in $10^{8}$. \label{fig:TwoPoreSep}}
\end{figure}

In Fig.~\ref{fig:TwoPoreSep} we show favorable comparisons between the numerical flux and the value of the series \eqref{ExactStriederResults}. Fig.~\ref{fig:TwoPoreSepA} shows the fluxes calculated from both methods and highlights the significant effect of interpore competition when pores are in close proximity. In Fig.~\ref{fig:TwoPoreSepB}, we observe the numerical method accurately resolves the $\bigoh (d^{-6})$ error term from the series solution \eqref{ExactStriederResults}. The method resolves errors to roughly one part in $10^{8}$.

\subsubsection{Example: Square and Hexagonal Pore arrangements}\label{Ex:SquareHexagon} 

\begin{figure}[htbp]
\centering
\subfigure[Square Pattern]{\includegraphics[width=0.475\textwidth]{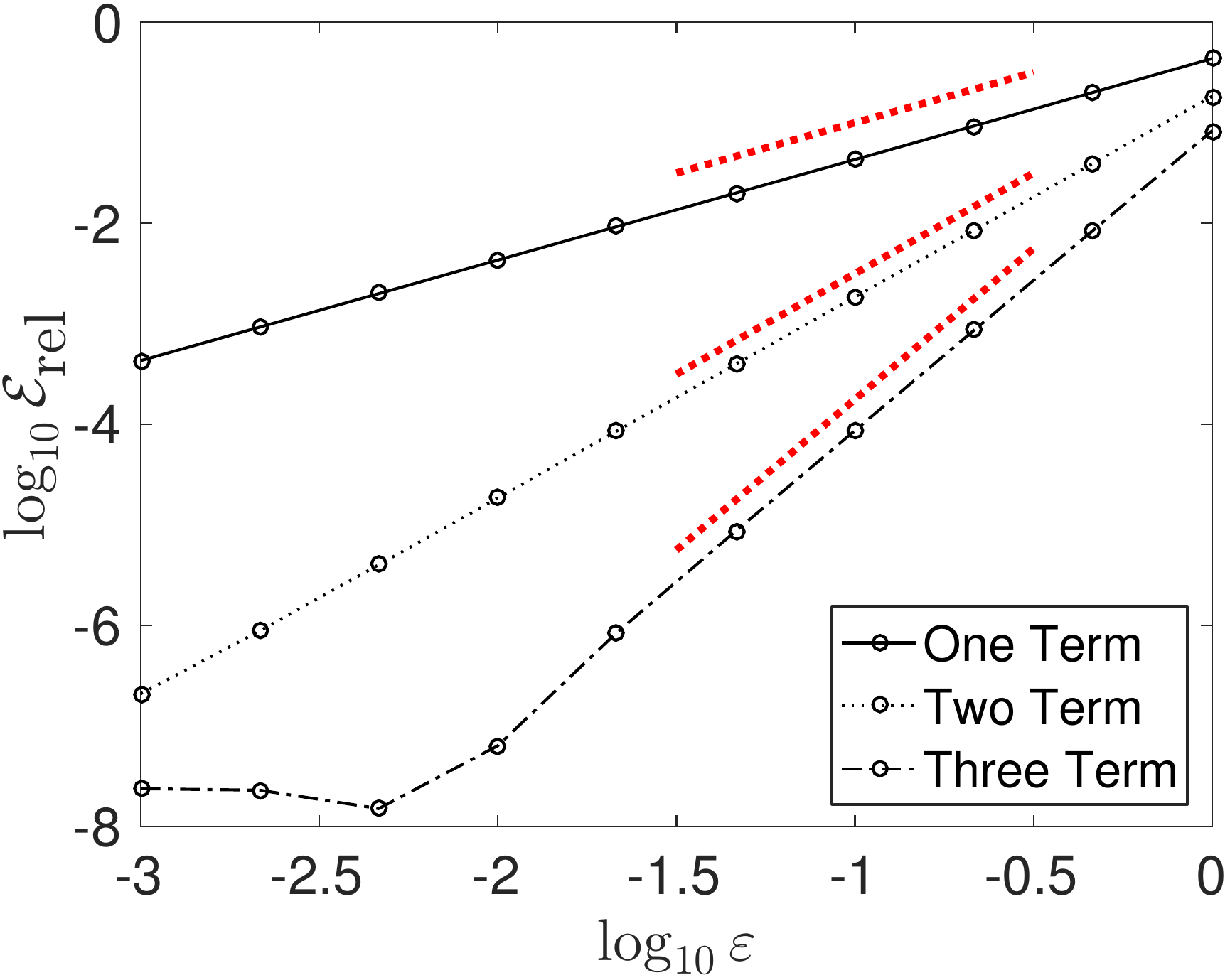} \label{fig:Square}}\quad
\subfigure[Hexagonal Pattern]{\includegraphics[width=0.475\textwidth]{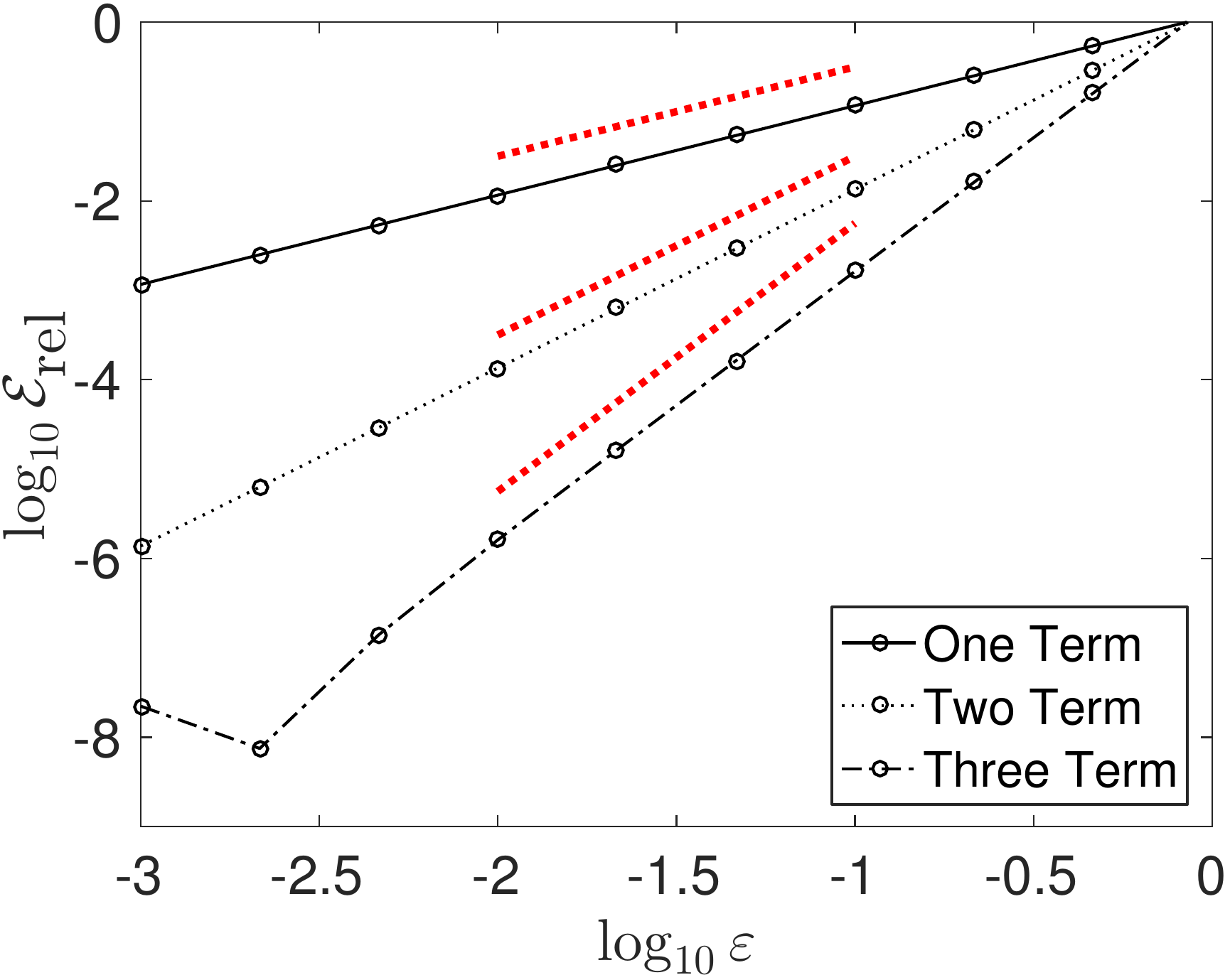} \label{fig:Hexagon}}
\caption{Results of example \ref{Ex:SquareHexagon}. Convergence of the numerical relative error with respect to the asymptotic approximation \eqref{mainResult} as the pore radius $\eps\to0$. The method is accurate to relative errors of around $10^{-8}$. Red dashed lines indicate lines of slope $1$, $2$ and $3$ corresponding to the error of the one, two and three term asymptotic approximations \eqref{mainResult}. \label{fig:SquareHexagon}}
\end{figure}

We verify the numerical method against asymptotic approximations in the limit of vanishing pore radius for a square and hexagonal planar pattern. In the square case the pore centers are $\bx = (\pm2,\pm2,0)$ while for the hexagonal case they are equally spaced on a ring of radius $2$. The pores have common radius $\eps$ which is varied and the relative error in the flux to the asymptotic prediction \eqref{mainResult}. Results in Fig.~\ref{fig:SquareHexagon} for $M=10$ polynomials show the numerical method is accurate to relative errors of around $10^{-8}$.

\begin{figure}[htbp]
\centering
\subfigure[Rescaled flux $J_s/(4\eps)$ against pore radius $\eps$.]{\includegraphics[width = 0.4725\textwidth]{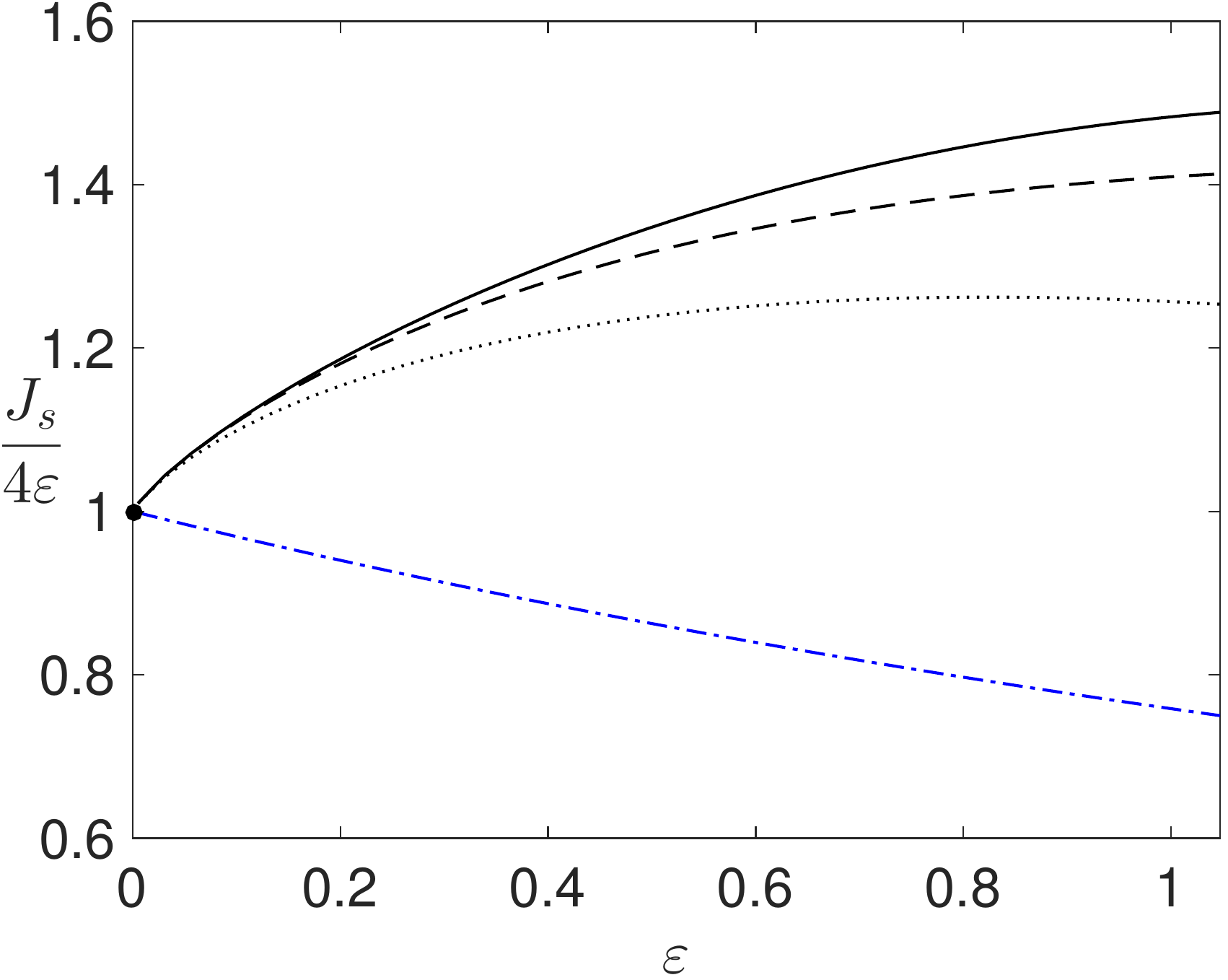}\label{fig:compareonepore_a}}\qquad
\subfigure[Relative errors $\mathcal{E}_{\text{rel}}$ on logarithmic scale.]{\includegraphics[width = 0.4725\textwidth]{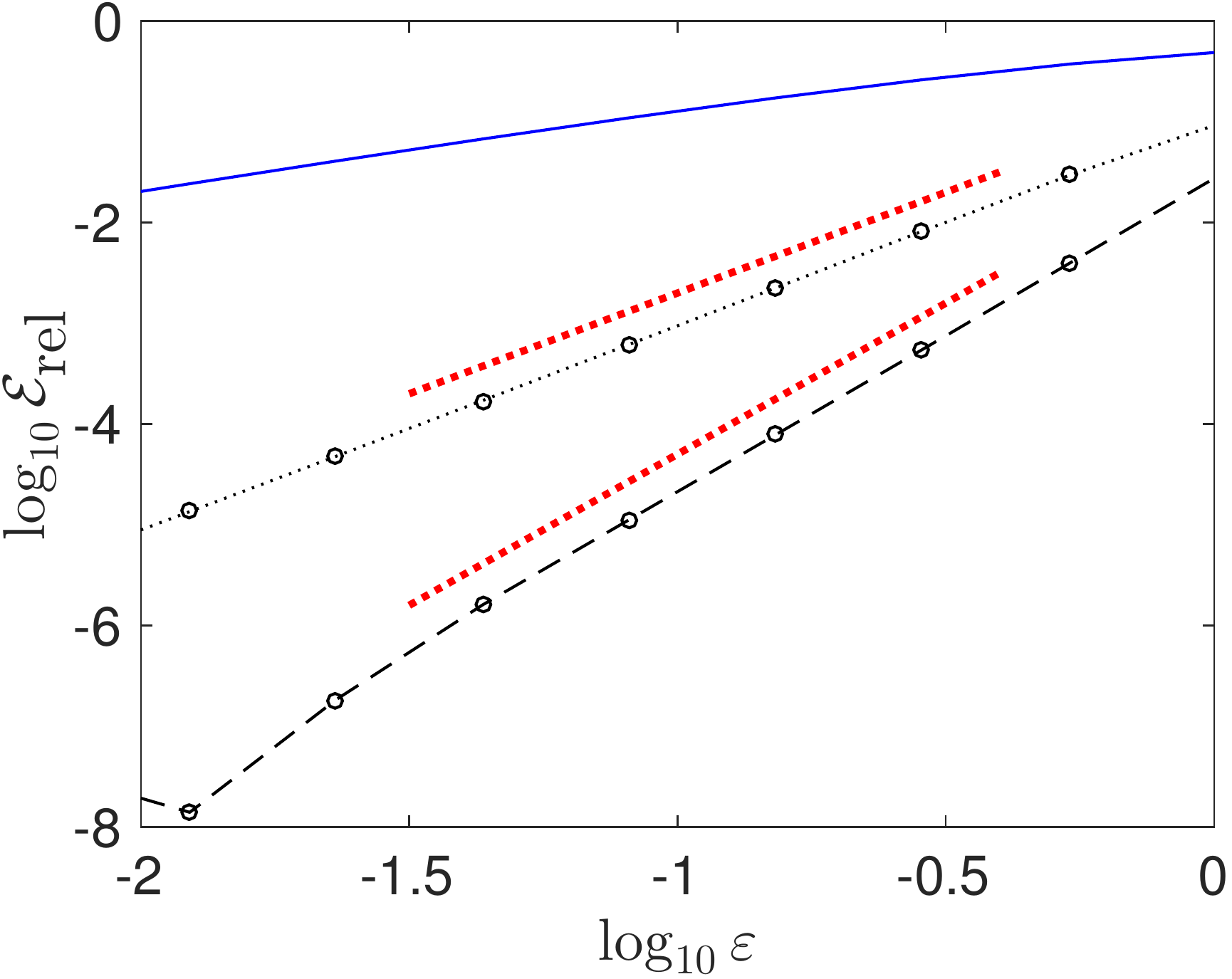}\label{fig:compareonepore_b}}
\parbox{0.92\textwidth}{\caption{Results for example \ref{sec:SPore} with a single pore on the sphere. Comparison of two and three term asymptotic predictions \eqref{CapOne} and numerics for a single pore of radius $\eps$ calculated with $M=20$ modes. Left panel: two term (black-dashed) expansion, three term (black-dotted) expansion, numerics (black-solid) and the Burg-Purcell (blue dot-dashes) from \eqref{resultBP}. Right panel: relative errors of asymptotic approximations for the capacitance as $\eps\to0$. Curves are the Berg-Purcell \eqref{resultBP} (solid blue), two term (dotted) and three term (dashed) asymptotic expansions from \eqref{CapOne}. Red lines of slope 2 (upper) and 3 (lower) confirm the expected order of the error. \label{fig:compareonepore}}}
\end{figure}

\subsection{Sphere Case}\label{resultsSphere}

In this section we consider the application of the numerical method to the spherical case. In \S\ref{sec:SPore}-\ref{Ex:PlantoicSpiral}, the numerical solution is validated using known closed form and asymptotic solutions, in the limit of small pore size and as the number of Zernike modes increases. Finally, in \S\ref{sec:Homogenization} we numerically validate a recently derived homogenized result which predicts the flux in terms of surface receptor density and typical pore size. Such results are crucial for use by experimentalists in real biological problems where the number of individual receptors is large and precise measurement of spatial locations impractical \cite{lag2,reing,WebsterCell2009}.

\subsubsection{Single pore}\label{sec:SPore}

For the single pore case, we verify the convergence of the spherical numerical method on single and multi pore cases. When $N=1$, a higher order approximation for the flux was derived in \cite{LWB2017} from a separable exact solution of \eqref{CapExt} 
\begin{equation}\label{CapOne}
J_s = 4D\eps \Big[1 + \frac{\eps}{\pi} \Big(
 \log 2\eps  -\frac32 \Big) - \frac{\eps^2}{\pi^2}\left( \frac{\pi^2 +21}{36}\right) + \mathcal{O}(\eps^3\log\eps)\, \Big]^{-1}\,, \qquad 
\mbox{as} \qquad \eps \to 0.
\end{equation}
 The results for the rescaled flux $J_s/(4\eps)$ in the case $N=1$, $D=1$ are shown in Fig.~\ref{fig:compareonepore_a} and demonstrate the validity of \eqref{CapOne}, even for moderately large pore radius. In Fig.~\ref{fig:compareonepore_b}, the numerical results give validation of the relative errors of the asymptotic formula as $\eps\to0$ and reveal that round off limits the smallest relative error obtainable to about one part in $10^8$. The blue curves in Fig.~\ref{fig:compareonepore} indicate the Berg-Purcell result \eqref{resultBP} which is significantly less accurate for $N=1$. 

\subsubsection{Example: Antipodal Pores}\label{Ex:TwoAntipodal}
\begin{figure}[htbp]
\centering
\subfigure[Antipodal pore configurations.]{\label{fig:antipodalA} \includegraphics[width=0.35\textwidth]{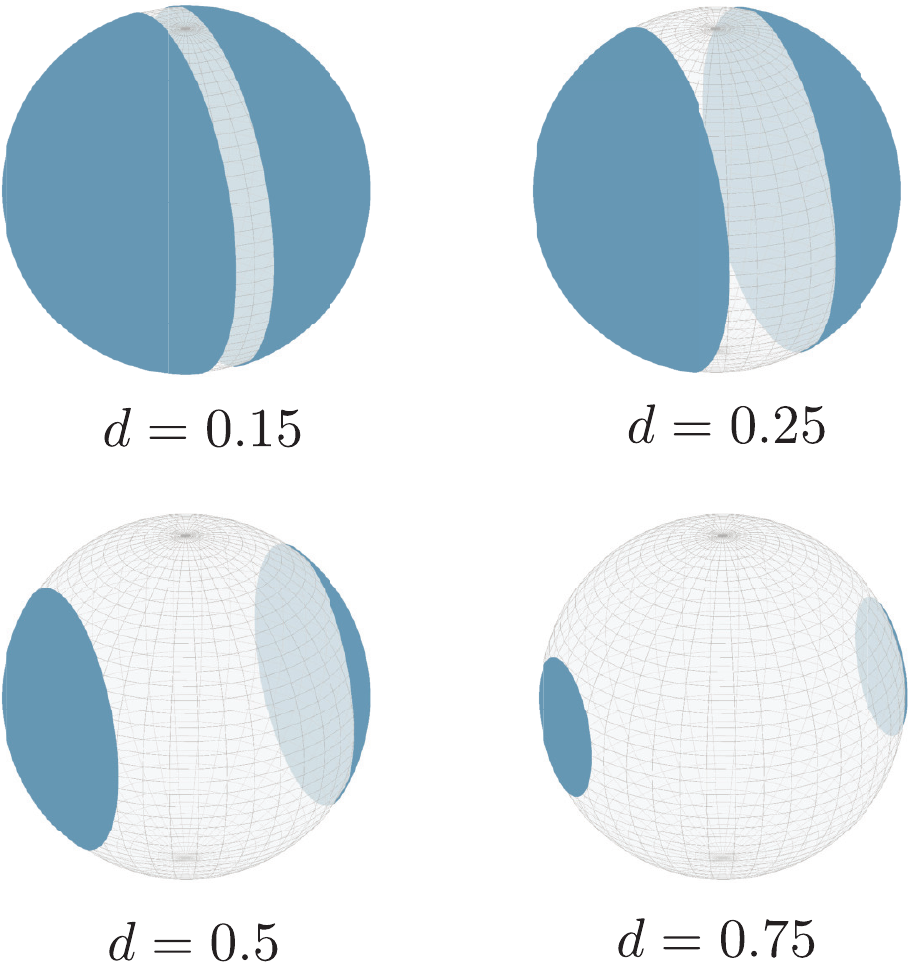}}\qquad
\subfigure[Relative error versus modes.]{\label{fig:antipodalB} \includegraphics[width=0.45\textwidth]{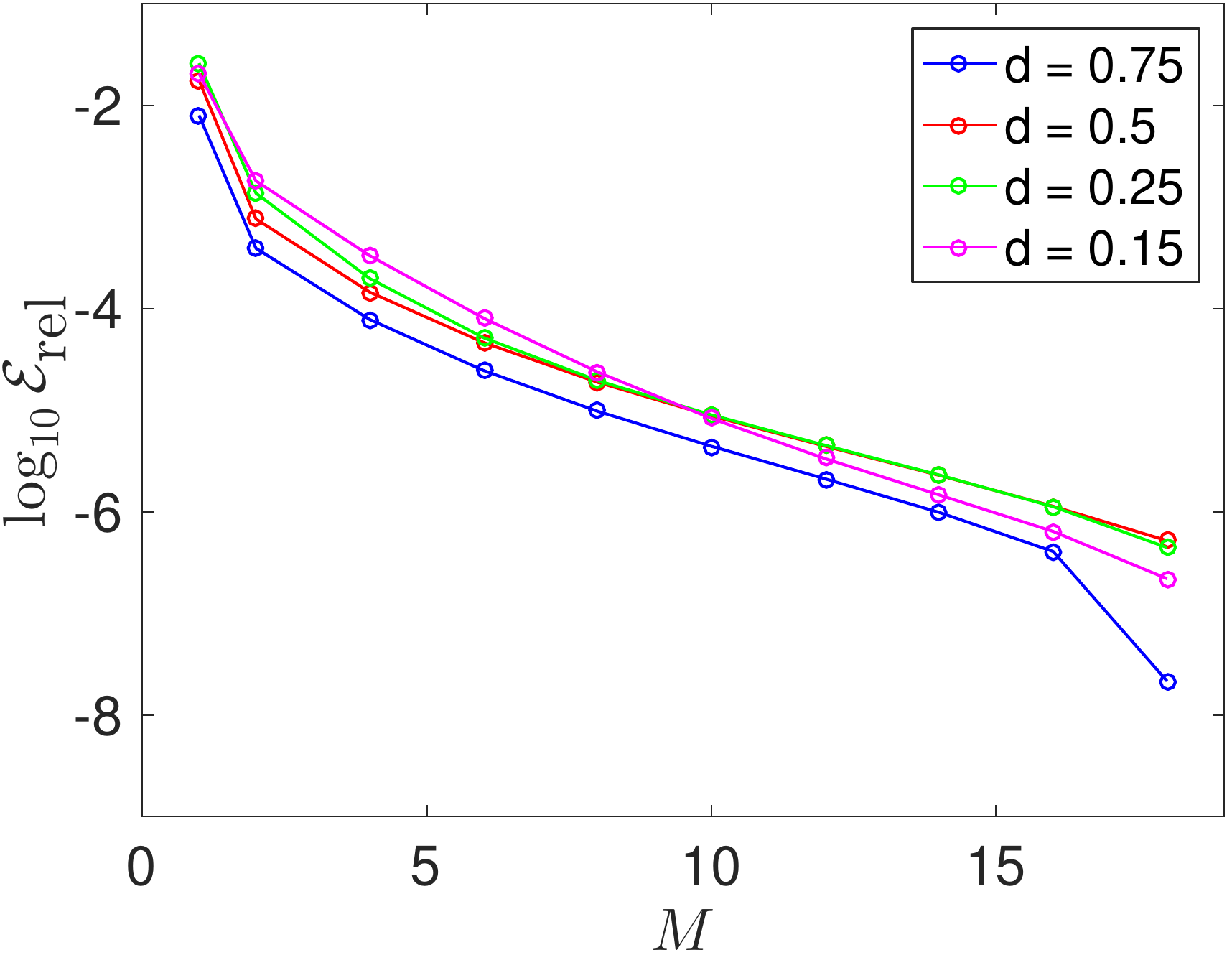}}
\parbox{0.92\textwidth}{\caption{Results of example \ref{Ex:TwoAntipodal} with $N=2$ antipodal pores of common radius $\alpha = \frac{\pi}{2}(1-d)$ for $d = 0.15, 0.25, 0.5$ and $d = 0.75$. As $d\to0$ and the interpore spacing decreases, additional modes are required to maintain numerical accuracy. Relative errors are calculated with respect to the \lq\lq true\rq\rq\ solution obtained from the spectral boundary element method evaluated with $M=20$ modes. \label{fig:antipodal}}}
\end{figure}

Here we consider $N=2$ pores in an antipodal position with common radius $\alpha=\frac{\pi}{2}(1-d)$ where $d$ is a separation parameter. For values $d= 0.15, 0.25, 0.5, 0.75$, we show convergence of the numerical flux as the number of modes increases. When $d\approx 1$ and the pore boundaries are well separated, fewer modes are necessary than for closely spaced pores $(d\approx 0)$. Results in Fig.~\ref{fig:antipodal}.

\begin{figure}[htbp]
\centering
\subfigure[Comparison with regular Platonic points.]{\label{fig:CompareMultiplePores_a} \includegraphics[width=0.45\textwidth]{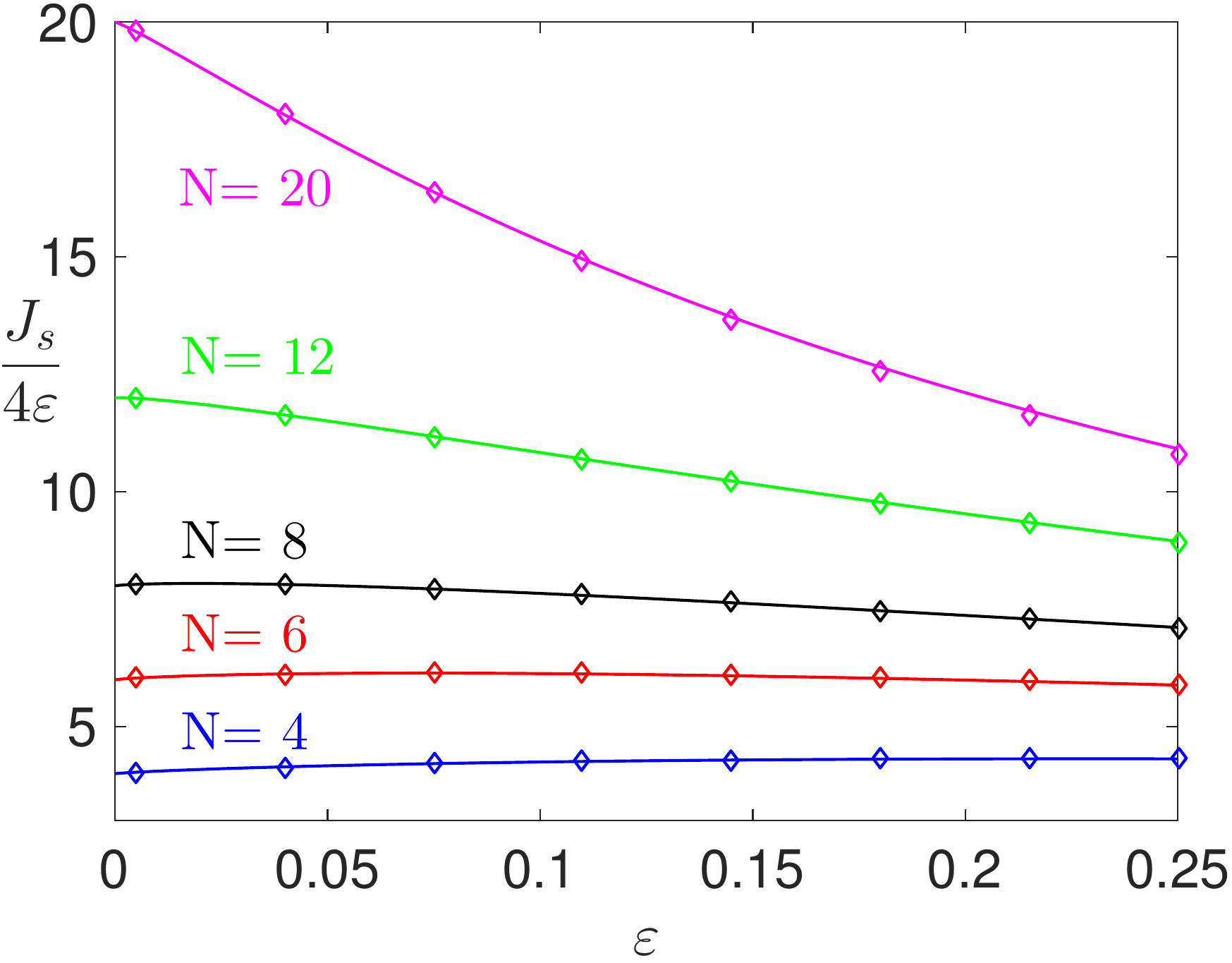}}\qquad
\subfigure[Comparison with Fibonacci lattice points.]{\label{fig:CompareMultiplePores_b} \includegraphics[width=0.45\textwidth]{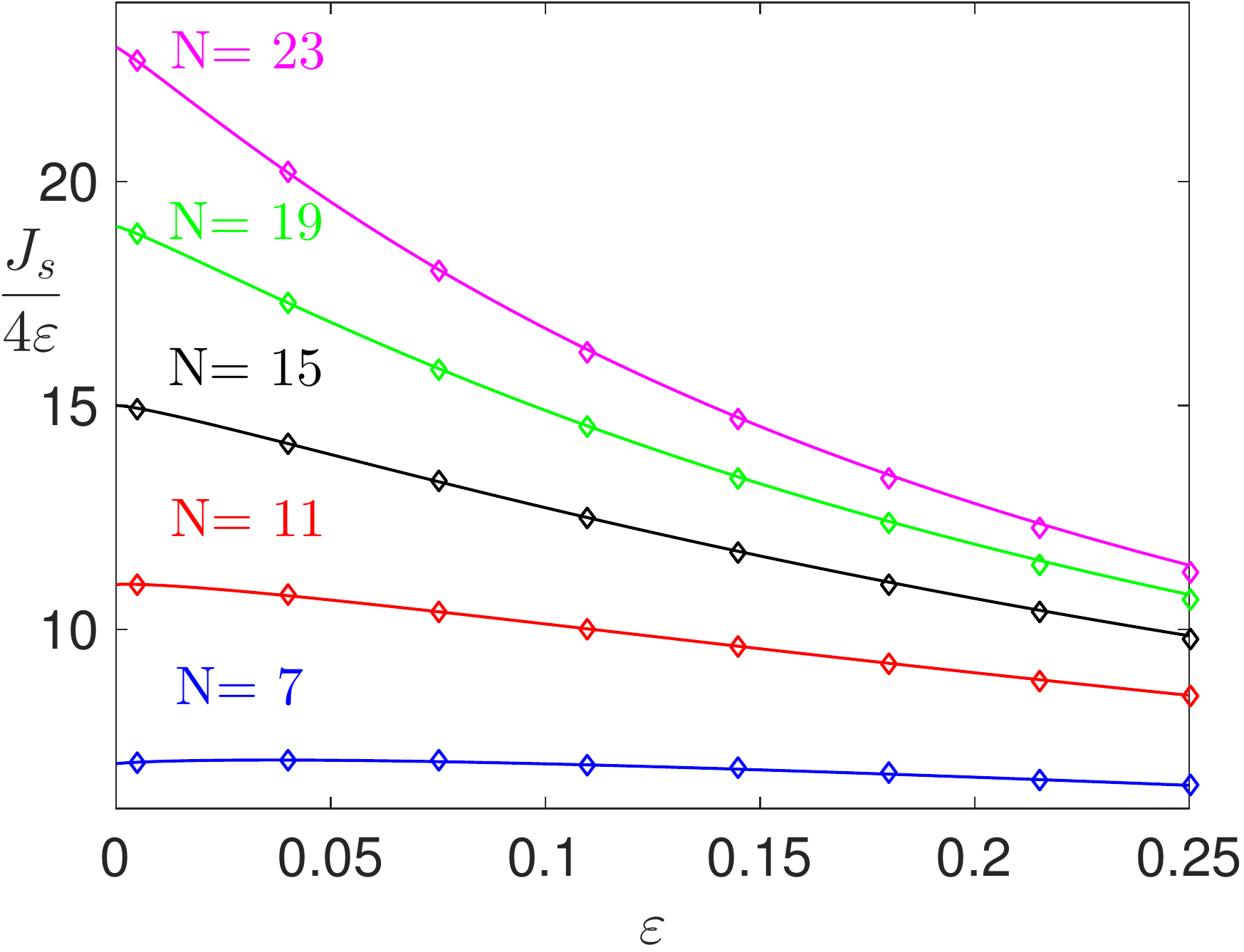}}
\parbox{0.92\textwidth}{\caption{Comparison of the rescaled flux $J_s/(4D\eps)$ as predicted by the asymptotic formula (solid lines) and full numerics (diamonds) using $M= 10$ modes and with pore locations given by vertices of the regular Platonic solids (left figure) and Fibonacci spirals (right figure). The Fibonacci spiral points \eqref{FibonacciPoints} generates an odd number of equispaced points on the sphere (see Fig.~\ref{Fig:Spiral}).\label{fig:CompareMultiplePores}}}
\end{figure}

\subsubsection{Example: Platonic Solids and Fibonacci Spirals}\label{Ex:PlantoicSpiral}

\begin{figure}[htbp]
\centering
\includegraphics[width=0.65\textwidth]{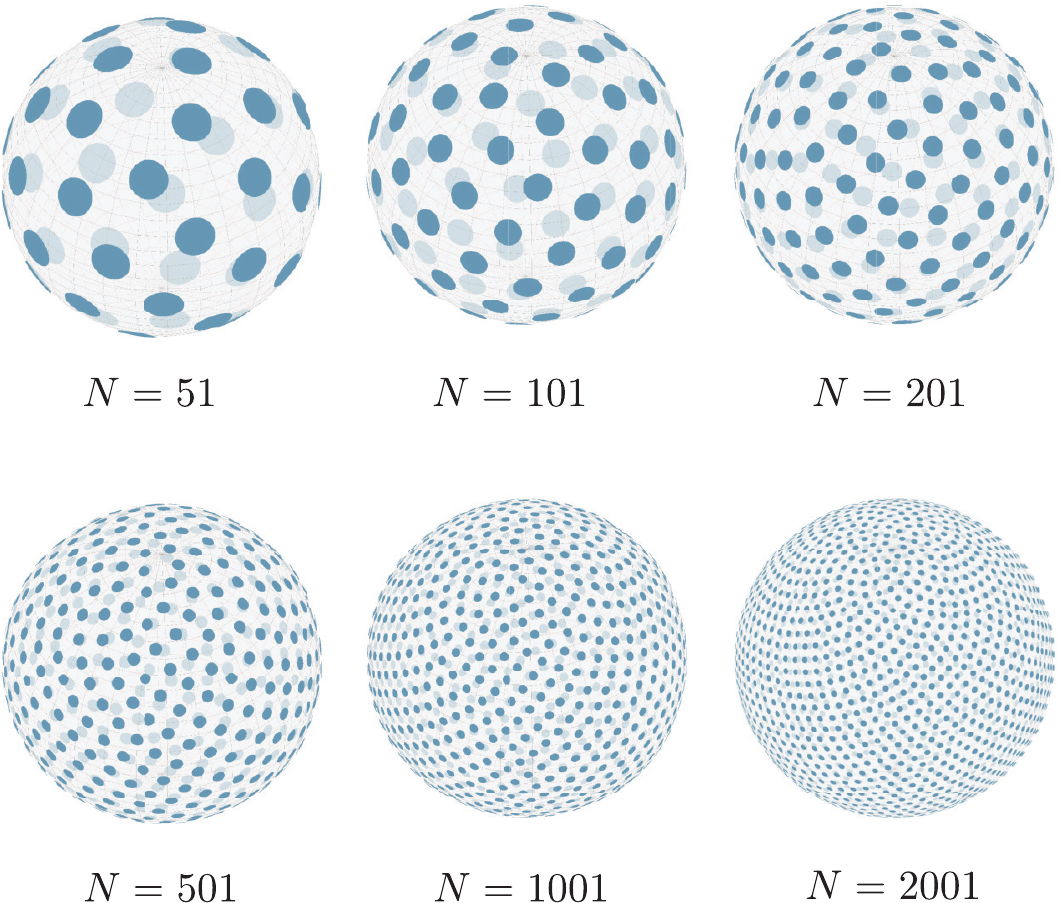}
\parbox{0.92\textwidth}{\caption{Homogeneous coverings of the sphere given by the Fibonacci spiral points \eqref{FibonacciPoints}. \label{Fig:Spiral} }}
\end{figure}

Here we verify the numerical method against the asymptotic approximation for the flux $J_s$ to multiple pores of common radius $\eps$ given by
\bsub\label{FluxSphereEqual}
\begin{equation}\label{FluxSphereEqual_a}
J_s = 4\eps D N\left[ 1 -  
 \frac{\eps}{\pi}\log 2\eps + \frac{\eps}{\pi} \Big(\frac{3}{2} -
  \frac{2}{N}\sum_{k\neq j} g_s(|\bx_j-\bx_k|) \Big) + \bigoh(\eps^2\log\eps) \right] \qquad \eps\to0,
\end{equation}
where the spherical pore interaction kernel $g_s(\mu)$ is given by
\begin{equation}\label{FluxSphereEqual_b}
g_s(\mu) = \frac{1}{\mu} + \frac{1}{2} \log\Big( \frac{\mu}{2+\mu} \Big), \qquad 0<\mu\leq2.
\end{equation}
\esub
The agreement between the numerical method and \eqref{FluxSphereEqual} is demonstrated for pores centered at the vertices of the regular Platonic solids in Fig.~\ref{fig:CompareMultiplePores_a} for $N\leq20$. We note that the vertices of the regular Platonic solids have many symmetries which can potentially obscure errors in the numerical method. It is therefore highly desirable to also benchmark the method against other distributions of spherical points. The equidistribution of a fixed number of points on the surface of a sphere is a long studied problem in approximation theory \cite{SK,GE,cerfon}. An easy to implement algorithm which produces a very homogeneously distributed set of points is the Fibonacci lattice \cite{Gonzalez2009,Swinbank2006}. Starting from an integer $k$, this algorithm produces $N= 2k+1$ points on the sphere with the $j^{\text{th}}$ point given in spherical coordinates by
\begin{equation}\label{FibonacciPoints}
\sin\theta_{j} = \frac{2j}{N}, \qquad \phi_j = \frac{2\pi j}{\Phi}, \qquad j = 1,\ldots, N,
\end{equation}
where $\Phi = 1 + \Phi^{-1} = (1+\sqrt{5})/2 \approx 1.618$ is the golden ratio. A few typical coverings arising from this algorithm are shown in Fig.~\ref{Fig:Spiral}. The accuracy of the numerical method with pores centered at Fibonacci vertices is demonstrated in Fig.~\ref{fig:CompareMultiplePores_b} by comparing to the asymptotic result \eqref{FluxSphereEqual}.

In Fig.~\ref{fig:CompareMultiplePores}, excellent agreement is seen in both cases for configurations up to $N=21$ pores with $M=10$ Zernike modes. As the common radius $\eps$ shrinks to zeros, we have from \eqref{FluxSphereEqual_a} that
\[
\lim_{\eps\to 0} \frac{J_s}{4D\eps}  = N,
\]
which agrees with the original Berg-Purcell result \eqref{resultBP} and is observed in each curve in Fig.~\ref{fig:CompareMultiplePores}.

\subsubsection{Homogenization}\label{sec:Homogenization}

In the cellular process of protein trafficking between the interior of the nucleus and the cytosol through Nuclear Pore Complexes (NPCs), the number of individual pores is approximately $N=2000$.  The nuclear radius is roughly $4$
microns and each NPC has an estimated radius of $25$ nanometres
(cf.~\cite{lag2,reing}). This implies that roughly $2\%$ of the
boundary of the nucleus is covered by pores.

It is experimentally impractical to accurately measure the 3D spatial location for each of the thousands of NPCs for use in \eqref{FluxSphereEqual}, however, the NPC density is comparatively simpler to obtain \cite{Maeshima10}. In the limit $\eps\to0$, $N\to \infty$, but with the absorbing surface area fraction $\sigma = (N\pi\eps^2)/(4\pi) = (N\eps^2)/4$ held fixed, a homogenized flux $J_h$ was derived in \cite{LWB2017} where
\begin{equation}\label{fluxHomogenized}
J_h = 4\pi D \left[ 1+ \frac{\pi\eps}{4\sigma} \left(  1- \frac{4}{\pi} \sqrt{\sigma} + \frac{\sigma}{\pi} \log( 4e^{-1} \sqrt{\sigma}) + \frac{\eps^2}{2\pi\sqrt{\sigma}} \right) \right]^{-1}.
\end{equation}
The homogenized formula \eqref{fluxHomogenized} was obtained from \eqref{FluxSphereEqual} assuming a uniform distribution of pores \cite{BS,RSZ1,KS,BBP,GE} with a combined absorbing surface area fraction satisfying $\sigma = \bigoh(-\eps^2\log\eps)$ as $\eps\to0$. To establish the accuracy of the formula \eqref{fluxHomogenized}, we simulate \eqref{CapExt} with up to $N=2001$ absorbing pores whose centers are the Fibonacci spiral points (cf.~Fig.~\ref{Fig:Spiral}). In table \eqref{tab:homogenized}, we find in the biological scenario highlighted above, the homogenized formula \eqref{fluxHomogenized} predicts the flux to the target to a relative error of approximately $0.34\%$.

\begin{table}[htbp]
\centering
\begin{tabular}{| c || c  c  c  c  c  c|}
\hline
$\sigma$ & N = 51 & N = 101 & N = 201 & N = 501 & N = 1001& N = 2001\\
\hline
$2\%$ & $1.02\%$ & $0.90\%$ & $0.76\%$ & $0.58\%$ & $0.37\%$ & $0.34\%$\\
$5\%$ & $1.29\%$  & $1.07\%$ & $0.87\%$ & $0.63\%$ & $0.48\%$ &$0.34\%$\\
$10\%$ & $1.42\%$ & $1.14\%$ &$0.90\%$ & $0.63\%$ & $0.47\%$ &$0.38\%$\\
$20\%$ &$1.43\%$  & $1.14\%$ & $0.89\%$ & $0.62\%$ & $0.46\%$ &$0.34\%$\\
\hline
\end{tabular}
\parbox{0.92\textwidth}{\caption{Percentage relative errors between the homogenized formula $J_h$ given in \eqref{fluxHomogenized} and the boundary element solution of \eqref{CapExt} calculated with $M=6$ modes for a range of pore surface area fractions $\sigma$. Pore centers given by the Fibonacci spiral points \eqref{FibonacciPoints} and shown in Fig.~\ref{Fig:Spiral}. \label{tab:homogenized}}}
\end{table}

\section{Discussion}\label{sec:discussion}

This paper has been concerned with the problem of determining the capture rate of three dimensional diffusing particles by absorbing surface pores. There are two main contributions. First, we have given explicit asymptotic expressions for the capture rate of diffusing particles by a finite collection of non overlapping absorbers arranged on either an infinite plane or the surface of a sphere. Second, we have introduced and validated a novel spectral boundary element method which provides a rapid and highly accurate numerical solution of this problem.

The analytical expressions for the capture rates give detailed information on the effect of clustering of receptor sites and the rate of capture of diffusing particles. Explicit results have previously only been obtained for the simplified scenario of one or two absorbers \cite{sneddon1966mixed,Strieder08,Strieder12,Holcman2006a}. 

The numerical method complements widely used particle based Monte Carlo methods. Its advantageous attributes are its high accuracy, quick runtime, and its recovery of a smooth solution to the underlying PDE \eqref{CapExt}. A limitation of the method is its explicit assumption of a circular pore geometry. Using this method, we have verified a recently derived homogenization result \eqref{fluxHomogenized} for the flux of particles to a spherical with numerous surface absorbers. In realistic biological scenarios in which $N\approx 2000$ pores occupying roughly $2\%$ of the surface area \cite{lag2,reing}, we find (cf.~table \ref{tab:homogenized}) that the homogenized theory predicts the flux to within a relative error of $0.34\%$.

There are many avenues of future investigation arising from this study. It would be highly desirable to obtain a homogenized theory directly from the the asymptotic result \eqref{AsyPlaneFlux} for pores centered at a variety of Bravais lattices \cite{Iron2014}. This would give a first principle derivation of the function form \eqref{KBerez} fitted in \cite{berez2012,berez2013,berez2014} by particle simulations. An extension of the spectral boundary element to periodic arrays of planar absorbers would be useful in accurately validating such homogenized theories. Finally, it is highly desirable to extend this work to sample the full distribution of capture times to a collection of small pores. This distribution describes the duration of a particle's search for a receptor and consequently sets the timescale of biophysical processes such as immune signaling. This problem is more challenging since it requires the solution of a parabolic equation in the exterior region, rather than the elliptic problem \eqref{CapExt}.
\section*{Acknowledgments}
A.~E.~Lindsay was supported by NSF grant DMS-1516753.  
A.~J.~Bernoff was supported by Simons Foundation grant 317319. The authors acknowledge useful conversations with M.~J.~Ward.

\bibliographystyle{siam}
\bibliography{newbib}

\end{document}